\documentclass{gtart}

\def\ifplaintex{\expandafter\ifx\csname documentclass\endcsname\relax}

\def\gtp{{\mathsurround=0pt\it $\cal G\mskip-2mu$eometry \&\ 
$\cal T\!\!$opology $\cal P\!$ublications}}  

\def\recd{{\small Received:\qua\receiveddate\ifx\reviseddate\relax
\else\qquad Revised:\qua\reviseddate\fi\par}} 

\def\lognumber#1{\def\thelognumber{#1}}
\def\volumenumber#1{\def\thevolumenumber{#1}}
\def\volumeyear#1{\def\thevolumeyear{#1}}
\def\papernumber#1{\def\thepapernumber{#1}}
\def\pagenumbers#1#2{\def\startpage{#1}\def\finishpage{#2}}
\def\published#1{\def\publishdate{#1}}

\def\received#1{\def\receiveddate{#1}}
\def\revised#1{\def\reviseddate{#1}}
\def\accepted#1{\def\accepteddate{#1}}
\def\asciititle#1{\def\theasciititle{#1}}
\def\covertitle#1{\def\thecovertitle{#1}}
\def\asciiauthors#1{\def\theasciiauthors{#1}}
\def\asciiaddress#1{\def\theasciiaddress{#1}}

\def\coverauthors#1{\def\thecoverauthors{#1}}
\long\def\asciiabstract#1{\long\def\theasciiabstract{#1}}

\let\\\par\let\thelognumber\relax\let\thevolumenumber\relax
\let\thepapernumber\relax\let\thevolumeyear\relax\let\startpage\relax
\let\finishpage\relax\let\publishdate\relax\let\receiveddate\relax
\let\reviseddate\relax\let\accepteddate\relax\let\theasciititle\relax
\let\thecovertitle\relax\let\theasciiauthors\relax\let\theasciiaddress\relax
\let\theasciiabstract\relax

\let\thecoverauthors\relax\let\theasciiemail\relax

\ifplaintex
\font\logobig=cmssbx10 scaled 3836
\font\logomed=cmssbx10 scaled 2557
\else
\font\logobig=cmssbx10 scaled 4200
\font\logomed=cmssbx10 scaled 2800
\fi

\long\def\makeagttitle{   
\count0=\startpage
\agt\hfill      
\hbox to 45truept{\vbox to 0pt{\vglue -13truept{\logomed A\kern -.37em{\logobig 
T}\kern -.38em G}\vss}\hss}
\break
{\small Volume \thevolumenumber\ (\thevolumeyear)
\startpage--\finishpage\nl
Published: \publishdate}

\vglue .25truein

{\parskip=0pt\leftskip 0pt plus
1fil\def\\{\par\smallskip}{\Large\bf\thetitle}\par\medskip} \vglue
0.05truein

{\parskip=0pt\leftskip 0pt plus 1fil\def\\{\par}{\sc\theauthors}
\par\medskip}%
 
\vglue 0.03truein 

{\small\leftskip 25truept\rightskip 25truept{\bf Abstract}\stdspace\theabstract

{\bf AMS Classification}\stdspace\theprimaryclass
\ifx\thesecondaryclass\relax\else; \thesecondaryclass\fi\par
{\bf Keywords}\stdspace \thekeywords\par}\vglue 7truept

}   

\ifplaintex
\font\phead=cmsl9 scaled 950
\font\pnum=cmbx10 scaled 913
\font\pfoot=cmsl9 scaled 950
\headline{\vbox to 0pt{\vskip -4.5mm\line{\small\phead\ifnum
\count0=\startpage ISSN 1472-2739 (on-line) 1472-2747 (printed)
\hfill {\pnum\folio}\else\ifodd\count0\def\\{ }%
\ifx\theshorttitle\relax\thetitle\else\theshorttitle\fi\hfill{\pnum\folio}
\else\def\\{ and }{\pnum\folio}\hfill\ifx\theshortauthors\relax\theauthors
\else\theshortauthors\fi\fi\fi}\vss}}
\footline{\vbox to 0pt{\vglue 0mm\line{\small\pfoot\ifnum\count0=\startpage
\copyright\ \gtp\hfill\else
\agt, Volume \thevolumenumber\ (\thevolumeyear)\hfill\fi}\vss}}
\else
\font=cmsl9 scaled 1050
\font\lnum=cmbx10 
\font=cmsl9 scaled 1050
\makeatletter
\def\@oddhead{{\small\lhead\ifnum\count0=\startpage ISSN 1472-2739 
(on-line) 1472-2747 (printed)\hfill {\lnum\number\count0}\else\ifodd\count0
\def\\{ }\ifx\theshorttitle\relax \thetitle \else\theshorttitle\fi\hfill
{\lnum\number\count0}\else\def\\{ and }{\lnum\number\count0}
\hfill\ifx\theshortauthors\relax 
\theauthors\else\theshortauthors\fi\fi\fi}}\def\@evenhead{\@oddhead}
\def\@oddfoot{\small\lfoot\ifnum\count0=\startpage\copyright\ \gtp\hfill\else
\agt, Volume \thevolumenumber\ (\thevolumeyear)\hfill\fi}
\def\@evenfoot{\@oddfoot}
\makeatother
\fi
\let\maketitlepage\makeagttitle
\let\makeshorttitle\maketitlepage
\let\maketitle\maketitlepage

\newwrite\gtoutfile
\long\gdef\makeheadfile{  
{\def\\{, }\def\s{ }
\immediate\openout\gtoutfile head.xxx
\immediate\write\gtoutfile{To: math@arxiv.org}
\immediate\write\gtoutfile{Subject: put OR rep NNNNN:ppppp}
\immediate\write\gtoutfile{--text follows this line--}
\immediate\write\gtoutfile{Proxy-for: \ifx\theasciiauthors\relax
\theauthors\else\theasciiauthors\fi\s<\ifx\theasciiemail\relax\theemail\else\theasciiemail\fi>}
\immediate\write\gtoutfile{\noexpand\\}
\immediate\write\gtoutfile{Authors: \ifx\theasciiauthors\relax
\theauthors\else\theasciiauthors\fi}
{\def\\{ }\immediate\write\gtoutfile{Title: \ifx\theasciititle\relax
\thetitle\else\theasciititle\fi}}
\immediate\write\gtoutfile{Subj-class: GT or SG, GR etc}
\immediate\write\gtoutfile{MSC-class: \theprimaryclass\ifx\thesecondaryclass\relax\else, \thesecondaryclass\fi}
\immediate\write\gtoutfile{Journal-ref: Algebr. Geom. Topol. \thevolumenumber\s
(\thevolumeyear) \startpage-\finishpage}
\immediate\write\gtoutfile{Comments: Published by Algebraic and
Geometric Topology at}
\immediate\write\gtoutfile{\s\s\s  http://www.maths.warwick.ac.uk/agt/AGTVol\thevolumenumber/agt-\thevolumenumber-\thepapernumber.abs.html}
\immediate\write\gtoutfile{\noexpand\\}
\immediate\write\gtoutfile{}
\ifx\theasciiabstract\relax
\immediate\write\gtoutfile{\theabstract}\else
\immediate\write\gtoutfile{\theasciiabstract}\fi
\immediate\write\gtoutfile{}
\immediate\write\gtoutfile{\noexpand\\}
\immediate\write\gtoutfile{}
\immediate\closeout\gtoutfile}}  

\def\maketitlepage{\makeagttitle\makeheadfile}
\let\makeshorttitle\maketitlepage
\let\maketitle\maketitlepage

\lognumber{13}
\volumenumber{3}
\volumeyear{2003}
\papernumber{13}
\published{8 May 2003}
\pagenumbers{399}{433}
\received{29 November 2002}
\revised{3 May 2003}
\accepted{14 January 2003}

\usepackage{amsmath,amssymb}

\newtheorem{theo}{Th\'eor\`eme}[section]

\newtheorem{lem}[theo]{Lemme}

\newtheorem*{conj}{Conjecture de non-r\'ealisation forte \cite{Ku}}

\newtheorem{pro}[theo]{Proposition}

\newtheorem{cor}[theo]{Corollaire}

\theoremstyle{definition}

\newtheorem{pro-def}[theo]{Proposition - D\'efinition}
\newtheorem{defi}[theo]{D\'efinition}

\newtheorem*{rem}{Remarque}

\let\fin\endproof

\def \ker{\mathrm{Ker}}

\def \d{\mathrm{d}}

\def\n{\mathbb{N}}
\def \s{\mathrm{s}}
\def\u {\mathcal{U}}
\def\nil {\mathcal{N} il}

\def \B{\mathrm{B}}
\def \C{\mathrm{C}}
\def \F{\mathrm{F}}
\def \H{\mathrm{H}}
\def \N{\mathrm{N}}
\def \P{\mathrm{P}}
\def \R{\mathrm{R}}
\def \T{\mathrm{T}}
\def \Z{\mathbb{Z}}
\def \Hom{\mathrm{Hom}}

\def\Fun{\rm{F}(1)}

\def \wS{\widehat{\mathcal{S}}}

\def \lra {\longrightarrow}
\def \hra {\hookrightarrow}
\def \era {\twoheadrightarrow}

\def \colim{\mathrm{colim}}
\def \cS{\mathcal{S}}
\def \E{\mathrm{E}}
\def \hom{\mathbf{hom}}
\def \hompt{\mathbf{hom}_{\mathrm{pt}}}
\def \I{\mathcal{I}}
\def \Im{\mathrm{Im}}
\def \J{\mathcal{J}}
\def \proo{\mathrm{pro-}}
\def \pt{\mathrm{pt}}
\def \r{\rule{0ex}{1.5ex}}
\def \res{\mathrm{Res}^{\bullet}}
\def \ce{\mathrm{S}^{1}}
\def \Sim#1{\triangle[#1]}
\def \Sk{\mathrm{Sk}}
\def \Top{\mathcal{T}\!\!\mathit{op}}
\def \tor{\mathrm{Tor}}
\def \tot{\mathrm{Tot}}
\def \n{\mathbb{N}}

\def \fp{\mathbb{F}_p}

\def \dright{\def \arraystretch{.2} \begin{array}{c} \rightarrow \\ \rightarrow \end{array}}
\def \carre#1#2#3#4{\arraycolsep=.25em\begin{array}[t]{ccc} #1 & \to & #2 \\ \downarrow & & \downarrow \\ #3 & \to & #4 \end{array}}

\def \Sq {\mathrm{Sq}}

\def \ie{\textit{i.e.~}}

\def \underline#1{\underset{\widetilde{~}}{#1}}

\begin{document}

\title{Espaces profinis et probl\`emes
de r\'ealisabilit\'e}
\covertitle{Espaces profinis et probl\noexpand\`emes
de r\noexpand\'ealisabilit\noexpand\'e}
\authors{Fran\c{c}ois-Xavier Dehon\\G\'erald Gaudens}
\coverauthors{Fran\noexpand\c{c}ois-Xavier Dehon\\G\noexpand\'erald Gaudens}

\address{Laboratoire J.A. Dieudonn\'e, Universit\'e de Nice 
Sophia-Antipolis\\Parc Valrose - BP 2053 - 06101 Nice, France\\{\rm and}\\
Laboratoire Jean Leray (UMR 6629 du C.N.R.S.), Universit\'e de Nantes\\
BP 92208 - 44322 Nantes Cedex 3, France}
\asciiaddress{Laboratoire J.A. Dieudonne, Universite de Nice 
Sophia-Antipolis\\Parc Valrose - BP 2053 - 06101 Nice, France\\and\\
Laboratoire Jean Leray (UMR 6629 du C.N.R.S.), Universite de Nantes\\
BP 92208 - 44322 Nantes Cedex 3, Franc}

\asciititle{Espaces profinis et problemes de realisabilite}
\asciiauthors{Francois-Xavier Dehon and Gerald Gaudens}               

\email{dehon@math.unice.fr, gaudens@math.univ-nantes.fr}

\begin{abstract}
The mod $p$ cohomology of a space comes with an action of the Steenrod Algebra.
L. Schwartz \cite{Sc2} proved a conjecture due to N. Kuhn \cite{Ku} stating that if the mod $p$ cohomology of a space is in a finite stage of the Krull filtration of the category of unstable modules over the Steenrod algebra then it is locally finite.
Nevertheless his proof involves some finiteness hypotheses.
We show how one can remove those finiteness hypotheses by using the homotopy theory of profinite spaces introduced by F. Morel \cite{Mo}, thus obtaining a complete proof of the conjecture.
For that purpose we build the Eilenberg-Moore spectral sequence and show its convergence in the profinite setting.
\end{abstract}

\asciiabstract{The mod p cohomology of a space comes with an action of
the Steenrod Algebra.  L. Schwartz [A propos de la conjecture de non
realisation due a N. Kuhn, Invent. Math. 134, No 1, (1998) 211--227]
proved a conjecture due to N. Kuhn [On topologicaly realizing modules
over the Steenrod algebra, Annals of Mathematics, 141 (1995) 321--347]
stating that if the mod $p$ cohomology of a space is in a finite stage
of the Krull filtration of the category of unstable modules over the
Steenrod algebra then it is locally finite.  Nevertheless his proof
involves some finiteness hypotheses.  We show how one can remove those
finiteness hypotheses by using the homotopy theory of profinite spaces
introduced by F. Morel [Ensembles profinis simpliciaux et
interpretation geometrique du foncteur T, Bull. Soc. Math. France,
124 (1996) 347--373], thus obtaining a complete proof of the
conjecture.  For that purpose we build the Eilenberg-Moore spectral
sequence and show its convergence in the profinite setting.}

\def \supprimer{
La cohomologie {\it modulo} $2$ d'un espace topologique est munie
d'une action de l'alg\`ebre de Steenrod. L. Schwartz a
d\'emontr\'e \cite{Sc2} une conjecture due \`a N. Kuhn \cite{Ku},
qui affirme que si la cohomologie {\it modulo} $2$ d'un espace est
dans un cran fini de la filtration de Krull de la cat\'egorie des
modules instables sur l'alg\`ebre de Steenrod, alors elle est
localement finie. Sa preuve fait cependant intervenir certaines
hypoth\`eses de finitude. Nous montrons dans cet article que la
th\'eorie homotopique des espaces profinis de F. Morel \cite{Mo}
alli\'ee \`a la m\'ethode de L. Schwartz m\`ene \`a une
d\'emonstration compl\`ete de cette conjecture. Pour cela, nous
construisons une suite spectrale d'Eilenberg-Moore et montrons sa
convergence dans le cadre profini.
}


\maketitlepage

\section{Introduction}

La cohomologie \`a coefficients dans un corps fini d'un espace topologique est naturellement munie d'une structure de module instable sur
l'alg\`ebre de Steenrod \cite{SE}. La question se pose de caract\'eriser, parmi les modules instables, ceux qui sont isomorphes \`a la
cohomologie d'un espace topologique. De tels modules instables sont dit {\it topologiquement r\'ealisables}.

D\'eterminer si un module instable est r\'ealisable est en g\'en\'eral difficile. Par exemple, le fameux r\'esultat de J.F.~Adams sur les \'el\'ements de l'homotopie des sph\`eres d'invariant de Hopf un \cite{Ad} se paraphrase en la non r\'ealisabilit\'e de certains modules finis.

Un crit\`ere de r\'ealisabilit\'e simple, utilis\'e dans la suite, repose sur la remarque suivante~: la cohomologie d'un espace est aussi munie d'une structure d'alg\`ebre, laquelle est compatible avec l'action de l'alg\`ebre de Steenrod \cite{SE, Sc1}. En particulier, l'\'el\'evation au carr\'e y est donn\'ee par une op\'eration de Steenrod. En d'autres termes, l'existence d'une structure d'alg\`ebre compatible fournit une condition n\'ecessaire \`a la r\'ealisabilit\'e d'un module instable.

Dans un article r\'ecent \cite{Ku}, N. Kuhn \'enonce, et d\'emontre dans certains cas, des conjectures sur la structure de module instable de la cohomologie
des espaces. Voici celle qu'on \'etudie ici.

\subsection{La conjecture de non-r\'ealisation forte due \`a N. Kuhn}

Soit $p$ un nombre premier. Notons $\u$ la cat\'egorie des modules instables sur l'alg\`ebre de Steenrod {\it modulo} $p$.
Soit $\u_n$ la sous-cat\'egorie pleine de $\u$ dont les objets sont ceux annul\'es par $\bar{\T} ^n$, le foncteur de Lannes r\'eduit it\'er\'e $n$ fois (voir la partie \ref{theorie:lannes}). La suite $(\u_n)$ co\"{\i}ncide avec la filtration de Gabriel-Krull de la cat\'egorie $\u$ (voir la partie \ref{filtrations}).

Le premier cran $\mathcal{U} _0$ de cette filtration est la sous-cat\'egorie pleine de $\mathcal {U}$ dont les objets sont les
modules \emph{localement finis}, \ie les modules $M$ tels que pour tout \'el\'ement $x$ de $M$, il existe un entier $\ell$ tel que toute op\'eration de degr\'e sup\'erieur \`a $\ell$ op\`ere trivialement sur $x$.

N. Kuhn a propos\'e la conjecture suivante \cite[strong realization conjecture, p.324]{Ku}~:
\begin{conj} Soit $X$ un espace topologique. Si la cohomologie \emph{modulo} $p$ de $X$ est dans $\u_n$ pour un certain entier $n$, alors elle est localement finie. \end{conj}

\subsection{L'approche de L. Schwartz}

L. Schwartz a d\'emontr\'e cette conjecture
pour $p=2$ sous certaines restrictions de
finitude. Pr\'ecis\'ement, il a
d\'emontr\'e \cite[Th\'eor\`eme 0.1]{Sc3}~:

\begin{theo}{\rm\cite{Sc3}}\qua
\label{theosc}
Soit $X$ un espace, et soit $M$ sa
cohomologie \emph{modulo} $2$. Supposons que $M$ est dans $\u_n$ pour un certain entier $n$, et que
pour tout entier $s$ le quotient $M_s/M_{s+1}$
a un nombre fini de g\'en\'erateurs sur
l'alg\`ebre de Steenrod. Alors $M$ est
localement fini, c'est-\`a-dire $M\in \u
_0$.
\end{theo}

Dans l'\'enonc\'e, $M_s$ d\'esigne le sous-module instable de $M$ form\'e des \'el\'ements $s$-nilpotents (voir la partie \ref{filtrations}).

Rappelons la trame de la d\'emonstration
de L. Schwartz. On proc\`ede par l'absurde
en supposant qu'il
existe un espace
$X$ dont la cohomologie \textit{modulo} $2$ est dans un cran
fini $\u _n$ de la filtration de Krull sans \^etre
localement finie. On aboutit \`a une
contradiction en deux \'etapes:

La premi\`ere \'etape est {\it la r\'eduction de N. Kuhn~:~}
On montre qu'il existe alors un
(nouvel) espace $X$ dont la cohomologie est dans le
premier cran $\u _1$ de la
filtration de Krull sans \^etre localement finie.

La seconde \'etape repose sur l'utilisation de la {\it suite
spectrale d'Eilenberg-Moore} associ\'ee \`a la fibration des chemins de $X$.
On pr\'ecise alors la structure
particuli\`ere des modules qui sont dans le premier cran de la
filtration de Krull sans \^etre localement finis pour construire pour un certain entier $d$ et
pour tout entier $i$ assez grand des classes $\alpha _{i,d}$ dans
la cohomologie de $X$, qui v\'erifient les deux propri\'et\'es
suivantes~:

\begin{itemize}

\item Le cup carr\'e de $\alpha _{i,d}$ est nul pour
tout entier $i$ assez grand,

\item La classe $\alpha _{i,d}$ est exactement $d$-nilpotente pour tout
entier $i$ assez grand.

\end{itemize}

Les propri\'et\'es particuli\`eres de ces classes ${\alpha}_{i,d}$ montrent qu'elles
induisent dans la cohomologie de l'espace de lacets $\Omega X$, \`a travers l'homomorphisme de coin, des classes
$\alpha _{i,d-1}$, qui
v\'erifient~:

\begin{itemize}

\item Le cup carr\'e de la classe
$\alpha _{i,d-1}$ est nul pour
tout entier $i$ assez grand,

\item La classe $\alpha _{i,d-1}$ est exactement
$(d-1)$-nilpotente.

\end{itemize}

Par une r\'ecurrence descendante, on arrive \`a la contradiction qu'il
existe des classes
$\alpha _{i,0}$ dans $\mathrm{H}^{*}{\Omega}^d X$ qui, pour $i$ assez grand, sont simultan\'ement
$0$-nilpotentes (\ie r\'eduites au sens de la structure d'alg\`ebre de la cohomologie) et de cup-carr\'e nul. Ceci contredit l'existence d'une
structure d'alg\`ebre compatible avec la structure de module sur l'alg\`ebre de Steenrod sur $\mathrm{H}^{*}{\Omega}^d
X$.

L'utilisation de
l'interpr\'etation g\'eom\'etrique du foncteur
$\mathrm{T}$ et de la suite spectrale
d'Eilenberg-Moore dans la cat\'egorie des
espaces impose les hypoth\`eses de
finitude faites dans le th\'eor\`eme
\ref{theosc}. Cet article expose le moyen de
contourner ces hypoth\`eses en utilisant la th\'eorie de
l'homotopie des espaces profinis.

\subsection{De l'utilisation des espaces profinis dans les probl\`emes de non
r\'ealisabilit\'e}
\label{kuhnprof:kuhntop}

Le premier point est qu'on conna\^{\i}t une interpr\'etation g\'eom\'etrique \emph{inconditionnelle} du foncteur de Lannes \cite{Mo}, le prix \`a payer \'etant d'\'echanger la cat\'egorie des espaces par celle des espaces profinis.

Le second point est que nous construisons dans cet article une suite spectrale d'Eilenberg-Moore dans le cadre des espaces profinis. Il faut pour cela donner un sens \`a la notion d'espace de lacets d'un espace profini point\'e. Le formalisme ad\'equat pour d\'evelopper ces deux points est l'alg\`ebre homotopique.

L'approche utilis\'ee par L. Schwartz s'applique alors {\it mutatis mutandi} et permet de montrer le th\'eor\`eme suivant.

\begin{theo}
\label{theo:dg}
Soit $X$ un espace profini. Si la
cohomologie \emph{modulo} $2$ continue de $X$
est dans un cran fini de la filtration de Krull, alors elle est localement
finie.
\end{theo}

On sait que la cohomologie d'un espace co\"{\i}ncide avec la cohomologie continue de son compl\'et\'e profini (voir la section \ref{s:profini}).
Par cons\'equent, en corollaire du
th\'eor\`eme \ref{theo:dg}, on \'etablit la
\emph{conjecture de non r\'ealisation
forte} pour $p=2$.

\begin{theo}
\label{dgs}
Soit $X$ un espace. Si la cohomologie \emph{modulo} $2$ de
$X$ est dans un cran fini de la filtration de Krull, alors elle est
localement finie.
\end{theo}

\subsection{Organisation de
l'article}
\label{plan}

L'article est organis\'e de la
mani\`ere suivante. Apr\`es avoir rappel\'e
les points de la th\'eorie homotopique des espaces
profinis que nous utilisons, nous
d\'etaillons les constructions des
sommes amalgam\'ees homotopiques
et produits fibr\'es homotopiques, et en particulier des espaces de lacets. Nous
construisons alors la suite spectrale d'Eilenberg-Moore d'un
produit fibr\'e homotopique et en
explicitons la structure  alg\'ebrique.
Suivent  quelques rappels sur la th\'eorie de Lannes dans
le cadre des espaces profinis. Dans la
derni\`ere section, on expose, suivant la
m\'ethode de L. Schwartz,  la d\'emonstration
du \emph{Th\'eor\`eme \ref{theo:dg}}, ce qui
\'etablit la {\it conjecture de non
r\'ealisation forte} (\emph{Th\'eor\`eme
\ref{dgs}}) pour $p=2$. Dans un
appendice, on d\'etaille certains points techniques.


La preuve du th\'eor\`eme
{\ref{theo:dg}} est,
tant dans sa ligne que dans son aspect
technique, tr\`es semblable \`a celle de L.
Schwartz {\cite{Sc3}}. N\'eanmoins,
dans un soucis de clart\'e, nous avons
donn\'e dans cet article une
d\'emonstration aussi compl\`ete que possible.
Une des
raisons qui motivent ce choix est qu'il nous faut nous soustraire des
hypoth\`eses de finitude faites dans
\cite{Sc3}. En effet, on y trouve
{\cite[p. 523]{Sc3}: {\it `On fera partout l'hypoth\`ese que les modules instables consid\'er\'es sont de dimension finie en chaque degr\'e.'}} De
plus, {\cite[lemme 1.12, p.
530]{Sc3}}
n'est pas vrai sans hypoth\`eses de
finitude, ce qui nous emp\^eche de nous
appuyer sur {\cite[propositions
1.9 et 1.10]{Sc3}}. Malheureusement, ces
propositions sont utilis\'ees en divers
points cruciaux de l'article, par
exemple dans la d\'emonstration de
{\cite[proposition 3.2]{Sc3}} ou encore
dans le calcul final {\cite[p. 342]{Sc3}}.


\rk{Remerciements}
Le second auteur remercie
chaleureusement le \emph{Centre de
Recerca Matem\`atica} de Barcelone ainsi que
le laboratoire Jean-Alexandre Dieudonn\'e
(Universit\'e de Nice) pour leur accueil durant la pr\'eparation de ce travail.

Nous remercions vivement Lionel Schwartz pour l'int\'er\^et qu'il a pris pour notre travail.

\section{Th\'eorie homotopique des espaces profinis (d'apr\`es \cite{Mo})}\label{s:profini}

\subsection{Espaces profinis}

Une cat\'egorie $\I$ est dite petite si la collection de ses objets forme un ensemble. Elle est dite filtrante si :
\begin{itemize}
\item[(i)] pour toute paire d'objets $i,j$ de $\I$, il existe un objet $k$ et des morphismes $k\to i$ et $k\to j$.
\item[(ii)] pour toute paire de morphismes $i\dright j$ entre deux objets $i,j$ il existe un objet $k$ et un morphisme $k\to i$ tel que les compos\'ees $k\to i\dright j$ soient \'egales.
\end{itemize}

Soit $\mathcal{C}$ une cat\'egorie. Un pro-objet de $\mathcal{C}$ est la donn\'ee d'une petite cat\'egorie filtrante $\I$ et d'un foncteur de $\I$ dans $\mathcal{C}$. On dit aussi un diagramme filtrant d'objets de $\mathcal{C}$.
Un morphisme entre deux pro-objets $X(-):\I\to\mathcal{C}$ et $Y(-):\J\to\mathcal{C}$ est un \'el\'ement de l'ensemble $\lim_{i}\colim_{j}\Hom(X(j),Y(i))$.
Les pro-objets et leurs morphismes forment une cat\'egorie not\'ee $\proo\mathcal{C}$. (Voir par exemple \cite[appendice]{AM}.)

Un ensemble profini est un espace topologique compact totalement discontinu.
La limite dans $\Top$ d'un diagramme filtrant d'ensembles finis (discrets) est un ensemble profini. Inversement tout ensemble profini appara\^{\i}t canoniquement comme la limite filtrante de ses quotients finis $X/R$, $R$ d\'ecrivant les relations d'\'equivalence ouvertes sur $X$.
(Voir par exemple \cite[2.8]{Do}.)

Notre cat\'egorie d'espaces est celle des ensembles simpliciaux, qu'on note $\cS$ (voir par exemple \cite[chap. VIII]{BK}).

On appelle espace profini un objet simplicial dans la cat\'egorie des ensembles profinis. Les espaces profinis et leurs morphismes forment une cat\'egorie not\'ee $\wS$.

Notons $\cS_{\mathrm{f}}$ la sous-cat\'egorie pleine de $\cS$ (et de $\wS$~!) form\'ee des ensembles finis simpliciaux.
Le foncteur $\wS\to\proo\cS_{\mathrm{f}}$ qui associe \`a un espace profini $X$ le diagramme filtrant de ses quotients finis simpliciaux $X/R$ est une \'equivalence de cat\'egories d'inverse le foncteur limite $\proo\cS_{\mathrm{f}}\to\wS$ (cf \cite[lemma 2.3]{Q2}).

\subsection{Th\'eorie homotopique}

Soit $X$ un espace profini. On d\'efinit la cohomologie {\it modulo} $p$ continue de $X$, not\'ee $\H^* X$, comme l'homologie du complexe des cocha\^{\i}nes continues de $X$ \`a valeurs dans le groupe discret $\Z/p$.
Cette cohomologie est aussi la colimite filtrante des cohomologies {\it modulo} $p$ ordinaires des quotients finis simpliciaux de $X$.
Elle h\'erite en particulier de la structure que poss\`ede
la cohomologie {\it modulo} $p$ ordinaire d'un espace :
c'est une alg\`ebre instable sur l'alg\`ebre de Steenrod
(voir par exemple \cite{Sc1}).

La cat\'egorie $\wS$ poss\`ede une structure de cat\'egorie de mod\`eles ferm\'ee avec pour \'equivalences faibles les morphismes induisant un isomorphisme en cohomologie {\it modulo} $p$ continue et pour cofibrations les monomorphismes \cite[th\'eor\`eme 1]{Mo}.

On en d\'eduit imm\'ediatement une structure de cat\'egorie de mod\`eles sur la cat\'egorie point\'ee $\wS_{\pt}$ form\'ee des espaces profinis munis d'un point base.
Observons que l'oubli du point base $\wS_{\pt}\to\wS$ admet un adjoint \`a gauche : le foncteur qui associe \`a un espace profini $X$ la r\'eunion disjointe de $X$ et d'un point base, qu'on note $X_{+}$.

\subsection{Compl\'etion profinie}

Le foncteur d'oubli de la
topologie $\wS\to\cS$ admet un adjoint
\`a gauche : le foncteur de {\it
compl\'etion profinie} qui associe \`a un ensemble simplicial $X$ la limite dans $\wS$ de ses quotients finis simpliciaux $X/R$.
Par construction, les morphismes $X\to X/R$ induisent un isomorphisme de la cohomologie {\it modulo} $p$ continue du compl\'et\'e profini de $X$ dans la cohomologie {\it modulo} $p$ ordinaire de $X$.

\subsection{Pro-$p$-compl\'etion et r\'esolution fibrante}

Soit $X$ un ensemble simplicial. On note $\res(X)$ sa $\mathbb{F}_{p}$-r\'esolution cosimpliciale \cite[part I]{BK}, $\tot_{s}\res(X)$ le $s$-i\`eme espace total partiel associ\'e, et $\P^{t}\tot_{s}\res(X)$ la $t$-i\`eme troncature de Postnikov de l'espace $\tot_{s}\res(X)$.
Alors :
\begin{itemize}
\item[--] L'espace $\tot_{0}\res(X)$ est un $\mathbb{F}_{p}$-espace affine simplicial et les morphismes $\tot_{s+1}\res(X)\to\tot_{s}\res(X)$ sont des fibration principales sous l'action de $\mathbb{F}_{p}$-espaces vectoriels simpliciaux.
\item[--] Lorsque $X$ est un ensemble fini simplicial, les espaces $\tot_{s}\res(X)$ le sont \'egalement.
\item[--] Lorsque $X$ est connexe, les morphismes naturels $X\to\tot_{s}\res(X)$ induisent un isomorphisme $\colim_{s}\H^{*}\tot_{s}\res(X) \to \H^{*}X$ \cite[chap. III, \S 6]{BK}.
\end{itemize}

\smallskip On d\'eduit des deux premiers points que lorsque $X$ est un ensemble fini simplicial, les espaces $\P^{t}\tot_{s}\res(X)$ sont des espaces $p$-finis ; c'est \`a dire ils sont finis simpliciaux, fibrants dans $\cS$, ont un nombre fini de composantes connexes et les groupes d'homotopie de chacune d'entre elles sont des $p$-groupes finis, triviaux sauf pour un nombre fini d'entre eux.

On d\'eduit du troisi\`eme point que les morphismes naturels $X\to \P^{t}\tot_{s}\res(X)$ induisent un isomorphisme $\colim_{s,t} \H^{*}\P^{t}\tot_{s}\res(X)\to \H^{*}X$ lorsque $X$ est fini simplicial.

Supposons maintenant que $X$ est un espace \emph{profini}. On note $\tot_{s}X$ l'espace profini $\lim_{\R} \tot_{s} \res(X/R)$.
Alors $\tot_{0}X$ est un $\mathbb{F}_{p}$-espace affine profini simplicial et les morphismes $$\tot_{s+1}X\to\tot_{s}X$$ sont des fibrations principales sous l'action de $\mathbb{F}_{p}$-espaces vectoriels profinis simpliciaux.
En particulier les morphismes
$\tot_{0}X\to \pt$ et
$\tot_{s+1}X\to\tot_{s}X$ sont des
fibrations dans $\wS$ (voir \cite[2.1]{Mo}).

On note $\widehat{X}(-)$ le diagramme filtrant form\'e des espaces $p$-finis $\P^{t} \tot_{s} \res (X/R)$, $t$ et $s$ d\'ecrivant l'ensemble des entiers positifs et $R$ l'ensemble des relations d'\'equivalence simpliciales ouvertes de $X$, et $\R X$ la limite de $\widehat{X}(-)$ dans $\wS$.
Le morphisme canonique $X\to\R X$ est une \'equivalence faible faisant de $\R X$ une r\'esolution fibrante fonctorielle de $X$.
On appelle le diagramme $\widehat{X}(-)$
le pro-$p$-compl\'et\'e de $X$.

\subsection{Objets fonctionnels}

Soient $X$ et $Y$ deux espaces profinis point\'es. On note $X\vee Y$ leur somme et $X\wedge Y$ le quotient $(X\times Y)/(X\vee Y)$.
L'espace profini $X\wedge Y$ est fonctoriel en $X$ et $Y$.

Notons, pour $W$ un espace profini point\'e, $\bar{\H}\r^* W$ sa cohomologie {\it modulo} $p$ continue r\'eduite, noyau du morphisme $\H^* W\to \H^*\pt$.
On dispose d'un morphisme canonique
$\bar{\H}\r^* X\otimes
\bar{\H}\r^* Y\to \bar{\H}\r^*
(X\wedge Y)$ qui est un isomorphisme si $X$ et $Y$ sont des ensembles finis simpliciaux donc un isomorphisme en g\'en\'eral par commutation des colimites filtrantes aux produits tensoriels.

Il existe un ensemble simplicial point\'e $\hompt(X,Y)$ caract\'eris\'e par la bijection
$$\Hom_{\cS_{\pt}}(W,\hompt(X,Y))\simeq \Hom_{\wS_{\pt}}(X\wedge W,Y)$$
naturelle en $W\in(\cS_{\mathrm{f}})_{\pt}$, et faisant de $\wS_{\pt}$ une cat\'egorie de mod\`eles ferm\'ee simpliciale (cf \cite[II,\S 1,2]{Q}).
Le fait que pour toute cofibration $A\to B$ et pour toute fibration $X\to Y$ dans $\wS_{\pt}$ le morphisme
$$\hompt(B,X)\to \hompt(A,X)\times_{\hompt(A,Y)}\hompt(B,Y)$$
est une fibration de $\cS_{\pt}$ et une \'equivalence faible si de plus $A\to B$ ou $X\to Y$ est une \'equivalence faible dans $\wS_{\pt}$ vient par adjonction des deux propri\'et\'es suivantes :
\begin{itemize}
\item[--] Une \'equivalence faible dans $\cS_{\pt}$ entre ensembles finis simpliciaux point\'es est une \'equivalence faible dans $\wS_{\pt}$.
\item[--] Pour tout monomorphisme $W\to W'$ entre ensembles finis simpliciaux point\'es (plus g\'en\'eralement pour toute cofibration $W\to W'$ dans $\wS_{\pt}$) le morphisme
$$(A\wedge W')\cup_{A\wedge W} (B\wedge W) \to B\wedge W'$$
est une cofibration de $\wS_{\pt}$ qui est une \'equivalence faible si $W\to W'$ ou $A\to B$ est une \'equivalence faible.
\end{itemize}

Supposons maintenant $X$ fini simplicial et notons $\Sk_{n}X$ son $n$-i\`eme squelette. L'ensemble simplicial $\hompt(X,Y)$ est alors naturellement profini car limite filtrante des ensembles finis simpliciaux $\hom_{\pt}(\Sk_{n}X,Y/R)$, et comme tel l'adjoint \`a droite en $Y$ du foncteur $Z\mapsto X\wedge Z$, $\wS_{\pt}\to \wS_{\pt}$.

La m\^eme discussion vaut dans $\wS$ en rempla\c{c}ant $X\wedge Y$ par $X\times Y$.
Observons que si $Y$ est point\'e l'ensemble simplicial $\hom(X,Y)$ est naturellement point\'e par le morphisme constant $X\to\pt\to Y$ et s'identifie comme tel \`a l'ensemble simplicial point\'e $\hompt(X_{+},Y)$.

{\bf Cas particulier}\qua Notons $\ce$ le quotient du simplexe standard de dimension $1$ de $\cS$ par son $0$-squelette.
L'espace $\ce$ est un ensemble fini simplicial point\'e.
A tout espace profini point\'e $X$ on associe sa suspension $\Sigma X=\ce\wedge X$ et son espace de lacets $\Omega X=\hompt(\ce,X)$.

Nous terminons par les trois remarques suivantes :
\begin{enumerate}

\smallskip\item[--] Soient $Y$ un espace profini point\'e, $X$ un ensemble simplicial point\'e quelconque et notons $(X_{\alpha})$ le diagramme filtrant form\'e des sous-ensembles finis simpliciaux de $X$.
On d\'efinit alors l'espace profini fonctionnel $\hompt(X,Y)$ comme la limite des espaces profinis $\hompt(X_{\alpha},Y)$.
Il est l'adjoint \`a droite en $Y$ du foncteur $\wS_{\pt}\to \wS_{\pt}$, $Z\mapsto\mathrm{colim}_{\alpha}(X_{\alpha}\wedge Z)$.

\smallskip\item[--] Soient $X$ un ensemble simplicial point\'e ayant un nombre fini de simplexes non d\'eg\'en\'er\'es et $Y$ un espace $p$-fini point\'e, alors l'espace $\hompt(X,Y)$ est \'egalement $p$-fini.

\smallskip\item[--] Soient $Y\to Y'$ une \'equivalence faible entre espaces profinis point\'es fibrants et $X$ un ensemble simplicial (respectivement un espace profini) point\'e, alors le morphisme $\hompt(X,Y)\to \hompt(X,Y')$ est une \'equivalence faible entre espaces profinis (respectivement entre ensembles simpliciaux), ceci parce qu'une \'equivalence faible entre espaces profinis (point\'es) fibrants est une \'equivalence d'homotopie simpliciale \cite[1.4, Lemme 3]{Mo}.
(Nous rappelons ci-dessous la notion d'homotopie simpliciale.)
\end{enumerate}

\section{Sommes amalgam\'ees et produits fibr\'es homotopiques}

\subsection{Sommes amalgam\'ees homotopiques}

Notons $\Sim{s}$ le simplexe standard de dimension $s$ de $\cS$ ; c'est un ensemble fini simplicial non fibrant dans $\cS$ si $s$ est non nul, faiblement \'equivalent au point.

Les morphismes $\d^{0},\d^{1}:\Sim{0}\to\Sim{1}$ et $\s^{0}:\Sim{1}\to\Sim{0}$ induisent pour tout espace profini point\'e $X$ des morphismes $X\vee X\to X\wedge\Sim{1}_{+}$ et $X\wedge\Sim{1}_{+}\to X$ qui font de $X\wedge\Sim{1}_{+}$ un (bon) objet cylindre pour $X$ dans $\wS_{\pt}$, fonctoriel en $X$.
(Voir par exemple \cite{DwS} pour cette notion).

Notons $\{0\}$ et $\{1\}$ les images respectives de $\d^{1}$ et $\d^{0}$ : $\Sim{0}\to\Sim{1}$.
Tout morphisme $X\to Y$ entre espaces profinis point\'es s'\'ecrit comme la compos\'ee des morphismes \'evidents $X\wedge\{0\}_{+}\to(X\wedge\Sim{1}_{+})\cup_{X\wedge\{1\}_{+}}Y$ et $(X\wedge\Sim{1}_{+})\cup_{X\wedge\{1\}_{+}}Y\to Y$.

\begin{pro}\label{Pcof} Soit $X\to Y$ un morphisme entre espaces profinis point\'es ; alors :
\begin{itemize}
\item[--] Le morphisme $X\wedge\{0\}_{+}\to(X\wedge\Sim{1}_{+})\cup_{X\wedge\{1\}_{+}}Y$ est une cofibration de $\wS_{\pt}$.
\item[--] Le morphisme $(X\wedge\Sim{1}_{+})\cup_{X\wedge\{1\}_{+}}Y\to Y$ est une \'equivalence faible de $\wS_{\pt}$.
\end{itemize}
\end{pro}

Soient maintenant $X,Y,Z$ des espaces profinis point\'es et $X\to Y$ et $X\to Z$ des morphismes.
On d\'efinit la somme amalgam\'ee homotopique du diagramme $Z\leftarrow X\rightarrow Y$ comme la somme amalgam\'ee dans $\wS_{\pt}$ du diagramme
$$(X\wedge\Sim{1}_{+})\cup_{X\wedge\{1\}_{+}}Z \leftarrow X\wedge\{0\}_{+}\rightarrow (X\wedge\Sim{1}_{+})\cup_{X\wedge\{1\}_{+}}Y\ .$$
Le foncteur qui associe au diagramme $Z\leftarrow X\rightarrow Y$ sa somme amalgam\'ee homotopique s'interpr\`ete comme le d\'eriv\'e gauche du foncteur qui associe \`a ce m\^eme diagramme sa colimite dans $\wS_{\pt}$ (cf \cite[10.7]{DwS}).

Supposons enfin que $X\to Y$ est injective degr\'e par degr\'e (\ie est une cofibration) et notons $\Sigma$ la somme amalgam\'ee dans $\wS_{\pt}$ du diagramme $Z\leftarrow X\rightarrow Y$.
Comme dans le cas classique le diagramme obtenu au niveau des complexes de cocha\^{\i}nes continues
$$\carre{\C^{*}\Sigma}{\C^{*}Z}{\C^{*}Y}{\C^{*}X}$$
est cocart\'esien de sorte qu'on obtient la suite exacte de Mayer Vietoris en cohomologie {\it modulo} $p$ continue.

En utilisant la fonctorialit\'e de cette suite exacte longue et la proposition qui pr\'ec\`ede on obtient que le morphisme canonique de la somme amalgam\'ee homotopique du diagramme $Z\leftarrow X\rightarrow Y$ dans $\Sigma$ est une \'equivalence faible.

{\bf Cas particulier}\qua Soit $X$ un espace profini point\'e et prenons pour $Y$ et $Z$ l'espace profini point\'e \'egal au point.
Le morphisme $(X\wedge\Sim{1}_{+})/(X\wedge\{1\}_{+}) \to \pt$ induit une \'equivalence faible de la somme amalgam\'ee homotopique du diagramme $\pt\leftarrow X\rightarrow \pt$ dans la suspension $\Sigma X$ de $X$.
La suite exacte de Mayer Vietoris devient l'isomorphisme en cohomologie {\it modulo} $p$ continue r\'eduite
$$\bar{\H}\r^* \Sigma X \cong \Sigma\bar{\H}\r^* X\ .$$

\subsection{Produits fibr\'es homotopiques}

Dualement notons $X^{\Sim{1}}$ l'espace profini $\hom(\Sim{1},X)$. Les morphismes $X\to X^{\Sim{1}}$ et $X^{\Sim{1}}\to X\times X$ induits par
$\d^{0},\d^{1}:\Sim{0}\to\Sim{1}$ et $\s^{0}:\Sim{1}\to\Sim{0}$ font de $X^{\Sim{1}}$ un espace de chemins pour $X$ dans $\wS$, fonctoriel en $X$.
(On peut se convaincre de ce que le morphisme $X\to X^{\Sim{1}}$ est une \'equivalence faible en lisant la d\'emonstration de la proposition qui suit.)
Tout morphisme $X\to Y$ entre espaces profinis s'\'ecrit comme la compos\'ee des morphismes $X\to X\times_{Y^{\{0\}}}Y^{\Sim{1}}$ et $X\times_{Y^{\{0\}}}Y^{\Sim{1}}\to Y^{\{1\}}$.

La proposition suivante est classique pour une cat\'egorie de mod\`eles ferm\'ee simpliciale et cruciale pour notre propos :

\begin{pro} \label{E} Soit $X\to Y$ un morphisme entre espaces profinis ; alors :
\begin{itemize}
\item[--] Le morphisme $X\to X\times_{Y^{\{0\}}}Y^{\Sim{1}}$ est une \'equivalence faible de $\wS$.
\item[--] Le morphisme $X\times_{Y^{\{0\}}}Y^{\Sim{1}}\to Y^{\{1\}}$ est une fibration de $\wS$ si $Y$ est fibrant.
\end{itemize}
\end{pro}

\rk{D\'emonstration}
La d\'emonstration du premier point est analogue \`a celle de la proposition 5 de \cite[II,\S 2]{Q} :
Il existe une homotopie simpliciale de la compos\'ee $\d^{1}\s^{0}:\Sim{1}\to\Sim{1}$ \`a l'identit\'e de $\Sim{1}$ relativement \`a $\d^{1}$, c'est \`a dire un
morphisme $\Sim{1}\times\Sim{1}\to\Sim{1}$ dont la restriction \`a $\Sim{1}\times\{0\}$ co\"{\i}ncide avec le compos\'e $\Sim{1}\to\{0\}\to\Sim{1}$, dont la restriction \`a $\Sim{1}\times\{1\}$ co\"{\i}ncide avec $\mathrm{Id}_{\Sim{1}}$ et dont la restriction \`a $\{0\}\times\Sim{1}$ se factorise par la projection $\{0\}\times\Sim{1}\to \{0\}\times\{0\}$.
On en d\'eduit que le morphisme $X\times_{Y^{\{0\}}}Y^{\{0\}}\to X\times_{Y^{\{0\}}}Y^{\Sim{1}}$ est une \'equivalence d'homotopie.
Or le fait que $X\times\Sim{1}$ soit un objet cylindre pour $X$ implique qu'une \'equivalence d'homotopie est une \'equivalence faible.

Pour le deuxi\`eme point supposons $Y$ fibrant ; alors le morphisme $Y^{\Sim{1}}\to Y^{\{0,1\}}$ est une fibration dans $\wS$.
On en d\'eduit que le morphisme $X\times_{Y^{\{0\}}}Y^{\Sim{1}}\to X\times_{Y^{\{0\}}}Y^{\{0,1\}}\simeq X\times Y^{\{1\}}$ est une fibration.
Maintenant si $X$ est fibrant alors la projection $X\times Y\to Y$ est une fibration.
\fin

On d\'efinit la fibre homotopique d'un diagramme $X\rightarrow Y\leftarrow Z$ d'espaces profinis comme le produit fibr\'e dans $\wS$ du diagramme
$$\R X\times_{\R Y^{\{0\}}}(\R Y)^{\Sim{1}} \rightarrow \R Y^{\{1\}} \leftarrow \R Z\times_{\R Y^{\{0\}}}(\R Y)^{\Sim{1}}\ .$$
(o\`u $\R X$, $\R Y$ et $\R Z$ d\'esignent les r\'esolutions fibrantes de $X$, $Y$ et $Z$).
Il s'interpr\`ete comme le d\'eriv\'e droit en $X\rightarrow Y\leftarrow Z$ du produit fibr\'e dans $\wS$ \cite[10.12]{DwS}.

\begin{pro}\label{E2}
Soit $X\rightarrow Y\leftarrow Z$ un diagramme entre espaces profinis. Supposons que $X$, $Y$ et $Z$ sont fibrants et que le morphisme $Z\to Y$ est une fibration, alors le morphisme de $X\times_Y Z$ dans le produit fibr\'e homotopique de $X\rightarrow Y\leftarrow Z$ est une \'equivalence d'homotopie.
\end{pro}

La proposition est cons\'equence de la proposition~\ref{E} et du lemme suivant \cite[II, lemma 8.10]{GJ}.

\begin{lem} Soient $X\rightarrow Y\leftarrow Z$ et $X'\rightarrow Y'\leftarrow Z'$ deux diagrammes entre espaces profinis fibrants et $X\to X'$, $Y\to Y'$, $Z\to Z'$ des \'equivalences faibles commutant avec les morphismes des diagrammes.
On suppose que les morphismes $Z\to Y$ et $X'\to Y'$ sont des fibrations de $\wS$. Alors le morphisme $X\times_Y Z\to X'\times_{Y'} Z'$ est une \'equivalence faible.
\end{lem}

{\bf Cas particulier}\qua Soit $Y$ un espace profini point\'e fibrant et prenons pour $X$ et $Z$ l'espace profini point\'e \'egal au point.
Les propositions \ref{E} et \ref{E2} montrent que le morphisme $\pt\to \pt\times_{Y^{\{0\}}}Y^{\Sim{1}}$ induit une \'equivalence faible de l'espace de lacets $\Omega Y$ dans le produit fibr\'e homotopique du diagramme $\pt\rightarrow Y\leftarrow \pt$.

Notre strat\'egie pour \'etudier la cohomologie {\it modulo} $p$ continue d'un produit fibr\'e homotopique est de comparer celui-ci avec la limite d'un diagramme filtrant de produits fibr\'es d'espaces $p$-finis auquel on associera un diagramme de suites spectrales d'Eilenberg-Moore.

Pr\'ecisons : On associe d'abord au diagramme $X\rightarrow Y\leftarrow Z$ le diagramme des pro-$p$-compl\'et\'es $\widehat{X}(-)\rightarrow \widehat{Y}(-)\leftarrow \widehat{Z}(-)$.
On modifie ensuite les cat\'egories index de $\widehat{X}(-)$, $\widehat{Y}(-)$ et $\widehat{Z}(-)$ de sorte que les morphismes $\widehat{X}(-)\to \widehat{Y}(-)$ et $\widehat{Z}(-)\to \widehat{Y}(-)$ soient des transformations naturelles entre diagrammes filtrants d'espaces (cf \cite[appendice, \S 3]{AM}).
Nous en rappelons l'argument :

Soient $X(-)$ un pro-objet de cat\'egorie index $\mathcal{I}$, $\mathcal{J}$ une autre cat\'egorie petite et filtrante et $\phi:\mathcal{J}\to\mathcal{I}$ un foncteur.
Les identit\'es de $X(\phi(j))$, $j$ d\'ecrivant $\mathcal{J}$, induisent un morphisme de pro-objets $X(-)\to X\circ\phi(-)$.

\begin{lem}\label{reind} Soit $X(-)\to Y(-)$ un morphisme entre pro-objets. Notons $\mathcal{I}_{X}$ et $\mathcal{I}_{Y}$ leurs cat\'egories index respectives.
Il existe une petite cat\'egorie filtrante $\mathcal{J}$, des foncteurs $\phi:\mathcal{J}\to \mathcal{I}_{X}$ et $\psi:\mathcal{J}\to\mathcal{I}_{Y}$ et une transformation naturelle $X\circ\phi\to Y\circ\psi$ tels que les morphismes $X(-)\to X\circ\phi(-)$ et $Y(-)\to Y\circ\psi(-)$ soient des pro-isomorphismes et tels que le diagramme de pro-objets
$$\carre{X(-)}{Y(-)}{X\circ\phi(-)}{Y\circ\psi(-)}$$
commute.
\end{lem}

{\bf D\'emonstration}\qua
Soit $\mathcal{J}$ la cat\'egorie d\'efinie comme suit :

\begin{itemize}

\smallskip\item ses objets sont les triplets $(i,j,X(i)\to Y(j))$, $(i,j)$ d\'ecrivant $\mathcal{I}_{X}\times\mathcal{I}_{Y}$, tels que le morphisme $X(i)\to Y(j)$ repr\'esente le morphisme $X(-)\to Y(j)$,

\smallskip\item un morphisme d'un objet $(i,j,X(i)\to Y(j))$ dans un objet $(i',j',X(i')\to Y(j'))$ est un couple de morphismes $(i\to i',j\to j')$ tels que le diagramme
$$\arraycolsep=.25em\begin{array}[t]{ccc}
X(i) & \rightarrow & Y(j) \\
\downarrow & & \downarrow \\
X(i') & \rightarrow & Y(j')
\end{array}$$
commute.

\end{itemize}

On note $\phi$ et $\psi$ les foncteurs \'evidents $\mathcal{J}\to \mathcal{I}_{X}$ et $\mathcal{J}\to \mathcal{I}_{Y}$.
Alors $\mathcal{J}$ est une cat\'egorie petite et filtrante, les images de $\mathcal{J}$ dans $\mathcal{I}_{X}$ et $\mathcal{I}_{Y}$ sont cofinales de sorte que les morphismes $X(-)\to X\circ\phi(-)$ et $Y(-)\to Y\circ\psi(-)$ sont des pro-isomorphismes (cf \cite[Appendice, \S 3]{AM}) et on dispose d'une transformation naturelle $X\circ\phi\to Y\circ\psi$ satisfaisant aux conditions de l'\'enonc\'e.
\fin

Observons que, $\phi$ \'etant fix\'e, le $\mathcal{J}$-diagramme $X\circ\phi$ est fonctoriel en le $\mathcal{I}_{X}$-diagramme $X$.
En it\'erant le proc\'ed\'e on voit que si $X(-)\rightarrow Y(-)\leftarrow Z(-)$ est un diagramme entre pro-objets, il existe une cat\'egorie filtrante $\mathcal{J}$, des foncteurs $\phi:\mathcal{J}\to\mathcal{I}_{X}$, $\psi:\mathcal{J}\to\mathcal{I}_{Y}$ et $\theta:\mathcal{J}\to\mathcal{I}_{Z}$ induisant des pro-isomorphismes $X(-)\to X\circ\phi(-)$, etc. et il existe des transformations naturelles $X\circ\phi\to Y\circ\psi$ et $Z\circ\theta\to Y\circ\psi$ tels que le diagramme
$${\arraycolsep=.25em\begin{array}[t]{ccccc}
X(-) & \rightarrow & Y(-) & \leftarrow & Z(-) \\
\downarrow & & \downarrow & & \downarrow \\
X\circ\phi(-) & \rightarrow & Y\circ\psi(-) & \leftarrow & Z\circ\theta(-)
\end{array}}$$
soit commutatif.

Revenons au produit fibr\'e homotopique d'un diagramme $X\rightarrow Y\leftarrow Z$ d'espaces profinis.
Ce qui pr\'ec\`ede permet de conclure que le diagramme $\R X\rightarrow \R Y\leftarrow \R Z$ est limite filtrante de diagrammes $\widehat{X}(\phi(j))\rightarrow \widehat{Y}(\psi(j))\leftarrow\widehat{Z}(\theta(j))$ entre espaces $p$-finis.
Le produit fibr\'e du diagramme
$$\R X\times_{\R Y} {\R Y}^{\Sim{1}} \rightarrow \R Y \leftarrow \R Z\times_{\R Y} {\R Y}^{\Sim{1}}$$
est alors la limite filtrante des produits fibr\'es des diagrammes entre espaces $p$-finis
$$\widehat{X}(\phi(j))\times_{\widehat{Y}(\psi(j))} \widehat{Y}(\psi(j))^{\Sim{1}} \rightarrow \widehat{Y}(\psi(j)) \leftarrow \widehat{Z}(\theta(j))\times_{\widehat{Y}(\psi(j))} \widehat{Y}(\psi(j))^{\Sim{1}}\ .$$

Nous somme pr\^ets \`a mettre en place la suite spectrale d'Eilenberg-Moore.

\section{Suite spectrale d'Eilenberg-Moore et cohomologie continue des espaces de lacets} \label{EM}

\subsection{Construction}

Nous suivons l'approche g\'eom\'etrique de Rector \cite{R}, voir aussi \cite{Bou} pour une comparaison avec la construction classique.

Soit $X\rightarrow Y\leftarrow Z$ un diagramme entre espaces profinis ; on lui associe le diagramme cosimplicial d'espaces profinis donn\'e par la construction cobar g\'eom\'etrique, qu'on note $\B^{\bullet}(X\rightarrow Y\leftarrow Z)$ ou plus simplement $\B^{\bullet}$ lorsqu'il n'y a pas d'ambigu\"{\i}t\'e : $\B^{n}=X\times Y^{n}\times Z$ pour $n\geq 0$, et la coaugmentation $X\times_Y Z\to \B^{0}$ (cf \cite{R}).

Rappelons qu'on d\'efinit le complexe normalis\'e $\N_* A$ d'un objet simplicial $A_{\bullet}$ d'une cat\'egorie ab\'elienne comme le complexe ayant pour objet en degr\'e $n$ l'intersection des noyaux des faces $\d_{i}$, $1\leq i\leq n$, et pour diff\'erentielle le morphisme induit par $\d_{0}$.
Ce complexe est isomorphe au complexe ayant pour objet en degr\'e $n$ le quotient de $A_n$ par l'image des d\'eg\'en\'erescences et pour diff\'erentielle le morphisme induit par la somme altern\'ee des faces.
Le $n$-i\`eme objet d'homologie $\H_n \N A$ du complexe $\N_* A$ s'identifie au $n$-i\`eme objet d'homotopie $\pi_n A$ de l'objet simplicial $A_\bullet$ et est isomorphe au $n$-i\`eme objet d'homologie du complexe $(A_n,\sum (-1)^i\d_i)$.
(On a une \'equivalence d'homotopie canonique entre les complexes $\N_* A$ et $A_*$.)
Voir par exemple \cite[\S 22]{Ma}.

On d\'efinit les espaces profinis point\'es $\N^0\B=\B^0_+$ et $\N^n\B=\B^n/(\Im(\d^1)\cup\ldots\cup\Im(\d^{n-1}))$ pour $n\geq 1$.
Le morphisme $\d^0:\B^n\to\B^{n+1}$ induit un morphisme $\N^n\B\to \N^{n+1}\B$ telle que le compos\'e $\N^n\B\to \N^{n+1}\B\to \N^{n+2}\B$ est le morphisme constant.
Le complexe induit en cohomologie continue r\'eduite par la suite $\N^0\B\to \N^1\B\to\cdots$ s'identifie au normalis\'e $\N_*\H^*\B$ du groupe ab\'elien gradu\'e simplicial $\H\r^*\B^\bullet$.

On construit ensuite inductivement une suite d'espaces profinis point\'es $X^n$ comme suit :
On pose $X^{-1}=(X\times_Y Z)_+$. La compos\'ee $X^{-1}\to\N^0\B\to\N^1\B$ est le morphisme constant.
Supposons construit $X^{n-1}$ avec un morphisme $X^{n-1}\to\N^n\B$ telle que la compos\'ee avec $\N^n\B\to\N^{n+1}\B$ soit constante.
Notons $N^{\prime n}$ l'espace profini point\'e $(X^{n-1}\wedge\Sim{1}_{+})\cup_{X^{n-1}\wedge\{1\}_{+}}\N^n\B$.
Le morphisme $X^{n-1}\to\N^n\B$ s'\'ecrit comme la compos\'e de la cofibration $X^{n-1}\to N^{\prime n}$ avec l'\'equivalence faible $N^{\prime n}\to \N^n\B$.
On d\'efinit $X^n$ comme le quotient de $N^{\prime n}$ par $X^{n-1}$.
La compos\'ee $N^{\prime n}\to \N^n\B\to\N^{n+1}\B$ se factorise par $N^{\prime n}\to X^n$ et $X^n\to\N^{n+1}\B$ v\'erifie l'hypoth\`ese de r\'ecurrence.

On obtient ainsi une suite de cofibrations $X^{n-1}\to N^{\prime n}\to X^n$ entre espaces profinis point\'es et d'\'equivalences faibles $N^{\prime n}\to\N^n\B$.

Posons $A^{-s,*}=\bar{\H}\r^{*+s+1}X^s$ pour $-s\leq 1$, $A^{-s,*}=A^{1,*}=\bar{\H}\r^* (X\times_Y Z)_+$ pour $-s>1$ et $E_1^{-s,*}=\bar{\H}^{*+s}\N^sB$.

On a un triangle exact

{\small\unitlength=10pt
$$\begin{picture}(17,3)(0,0)
\put(2,3){$A^{-s+1,*}$}\put(6.5,3.2){\vector(1,0){4.2}}
\put(11.7,3){$A^{-s,*}$}
\put(7.6,0.5){\vector(-1,1){2}}\put(11.1,2.5){\vector(-1,-1){2}}
\put(7,-1){$E_1^{-s,*}$}
\end{picture}
$$}

o\`u $A^{-s,*}\to E_1^{-s,*}$ est un morphisme de bidegr\'e $(0,1)$,
donc une suite spectrale $(E_r,\d_r:E_r^{-s,t}\to E_r^{-s+r,t+1})$ v\'erifiant : (Cf les lemmes 5.6, 5.9 et le th\'eor\`eme 6.1 de \cite{Boa}.)
\begin{itemize}
\item[--] Le terme $E_1^s$ est nul pour tout $s>0$ donc la suite $(E_r^{-s})_r$ est une suite d'\'epimorphismes pour $r>s$. On pose $E_\infty^{-s}=\mathrm{colim}_rE_r^{-s}$.
\item[--] Le groupe ab\'elien gradu\'e $A^1\cong\H^* X\times_Y Z$ a une filtration naturelle $0\subset \F_0A^1\subset \F_{-1}A^1\subset\cdots\subset A^1$ d\'efinie par $\F_{-s}=\ker(A^1\to A^{-s})$.
On a un \'epimorphisme $E_\infty^{-s}\era \F_{-s}A^1/\F_{-s+1}A^1$ et en particulier un morphisme de coin $E_r^0\era \F_0A^1\subset A^1$.
\item[--] Les conditions suivantes sont \'equivalentes :
\begin{itemize}
\item[(1)] La filtration de $A^1$ est exhaustive et le morphisme $E_\infty^{-s}\to \F_{-s}A^1/\F_{-s+1}A^1$ est un isomorphisme pour tout $s$.
\item[(2)] Le groupe ab\'elien gradu\'e $\mathrm{colim}_s A^{-s}$ est nul.
\end{itemize}
\end{itemize}

On r\'eindexe la suite spectrale en posant $\E_r^{-s,t}=E_r^{-s,t+s}$ de sorte que le terme $\E_1$ est donn\'e par
$$\E_1^{-s,t}=\bar{\H}\r^t\N^s\B \cong \N_s\H^t\B$$
et la diff\'erentielle $\d_r$ est de bidegr\'e $(r,1-r)$.

On appelle cette suite spectrale la suite spectrale d'Eilenberg-Moore en cohomologie {\it modulo} $p$ continue associ\'ee au diagramme $X\rightarrow Y\leftarrow Z$.

\subsection{Convergence}
\label{convergence}

Le th\'eor\`eme 3.1 de \cite{Sh}, qui s'appuie sur le th\'eor\`eme de convergence
forte de Dwyer \cite{Dw}, entra\^{\i}ne la suivante~:

\begin{pro} Soit $X\rightarrow Y\leftarrow Z$ un diagramme entre espaces $p$-finis tel que $Z\to Y$ est une fibration, alors la suite spectrale d'Eilenberg-Moore associ\'ee converge fortement vers la cohomologie {\it modulo} $p$ du produit fibr\'e $X\times_Y Z$.
\end{pro}

Soit maintenant $X\rightarrow Y\leftarrow Z$ un diagramme entre espaces profinis quelconques.
Le morphisme de ce diagramme dans sa ``r\'esolution fibrante" $\R X\times_{\R Y} {\R Y}^{\Sim{1}} \rightarrow \R Y \leftarrow \R Z\times_{\R Y} {\R Y}^{\Sim{1}}$ induit un isomorphisme au niveau des termes $\E_2$ des suites spectrales d'Eilenberg-Moore.
La proposition qui pr\'ec\`ede, la discussion sur les produits fibr\'es et l'exactitude des colimites filtrantes entra\^{\i}nent :

\begin{pro} Soit $X\rightarrow Y\leftarrow Z$ un diagramme entre espaces profinis. La suite spectrale d'Eilenberg-Moore associ\'ee converge vers la cohomologie {\it modulo} $p$ continue de son produit fibr\'e homotopique.
\end{pro}

\textbf{Cas particulier}\qua Soit $X$ un espace profini point\'e fibrant. La suite spectrale d'Eilenberg-Moore associ\'ee au diagramme $\pt\rightarrow X\leftarrow \pt$ converge vers la cohomologie {\it modulo} $p$ continue de l'espace profini de lacets $\Omega X$.

\subsection{Structure de la suite spectrale d'Eilenberg-Moore}

Rappelons que pour $X$ et $Y$ deux espaces profinis on dispose d'un isomorphisme canonique $\H^*X\otimes\H^*Y\to\H^*X\times Y$.

L'objet gradu\'e simplicial $\H^*\B^\bullet$ s'identifie donc au produit tensoriel au dessus de $\H^*Y$ de $\H^*X$ avec la r\'esolution simplicial canonique de $\H^*Z$ comme $\H^*Y$-module.
Notons $\bar{\H}\r^*Y$ le conoyau de l'unit\'e $\fp\to\H^*Y$.
Le complexe normalis\'e $\N_*\H^*\B$ s'identifie au produit tensoriel au dessus de $\H^*Y$ de $\H^*X$ avec la r\'esolution bar r\'eduite de $\H^*Z$ comme $\H^*Y$-module
(cf \cite[Chap. X, \S 2]{Mac}).
Autrement dit on a un isomorphisme canonique
$$\E_1^{-s,*}\cong \H^*X\otimes(\bar{\H}\r^*)^{\otimes s}\otimes\H^*Z\ .$$
Le terme $\E_2^{-s,*}$ s'identifie \`a l'objet $\tor_s^{\H^*Y}(\H^*X,\H^*Z)$ et on a un morphisme de coin
$$\E_2^{0,*}\cong \H^*X\otimes_{\H^*Y}\H^*Z\era \E_\infty^{0,*}\to\F_0\H^*X\times_Y Z\subset \H^*X\times_Y Z\ .$$

\textbf{Cas particulier}\qua Supposons $X=Z=\pt$ et $Y$ fibrant.
Le morphisme $X\to Y$ munit $Y$ d'un point base et on a un isomorphisme canonique $\bar{\H}\r^*Y\cong \ker(\H^*Y\to\H^*\pt)$.
Le terme $\E_1^{-s,*}$ s'identifie \`a $(\bar{\H}\r^*Y)^{\otimes s}$.
La diff\'erentielle $\d_1:\E_1^{-1,*}\to\E_1^{0,*}=\fp$ est nulle.
La diff\'erentielle $\d_1:\E_1^{-2,*}\to\E_1^{-1,*}$ est donn\'ee par le cup produit restreint \`a $\bar{\H}\r^*Y$.
La filtration de la cohomologie continue de $\Omega Y$ induit une filtration de la cohomologie continue r\'eduite de $\Omega Y$ et on dispose de morphismes de coin
$$\fp=\E_1^{0,*}=\E_\infty^{0,*}\cong\F_0\H^*\Omega Y=\fp$$
et
$$\bar{\H}\r^*Y\cong\E_1^{-1,*}\era\E_\infty^{-1,*}\cong\Sigma\F_{-1}\bar{\H}\r^*\Omega Y \ .$$

Venons en \`a la structure de module sur l'alg\`ebre de Steenrod.

Le triangle exact 
{\small\unitlength=10pt
$$\begin{picture}(17,3)(0,1)
\put(3,3){$A^{-s+1,*}$}\put(7.1,3.2){\vector(1,0){3}}
\put(11.1,3){$A^{-s,*}$}
\put(8,1){\vector(-1,1){1.5}}\put(10.5,2.5){\vector(-1,-1){1.5}}
\put(7.5,-.4){$E_1^{-s,*}$}
\end{picture}$$}

dont est issue la suite spectrale est un triangle exact de modules instables sur l'alg\`ebre de Steenrod.
Les groupes ab\'eliens gradu\'es $\E_r^{-s,*}$ sont donc des modules instables pour tous $r,s$ de fa\c{c}on compatible avec l'aboutissement et la diff\'erentielle $\d_r:\E_r^{-s,*}\to\Sigma^{r-1}\E_r^{-s+r,*}$ est un morphisme de modules instables.

Donnons nous maintenant deux diagrammes $X\rightarrow Y\leftarrow Z$ et $X'\rightarrow Y'\leftarrow Z'$ entre espaces profinis et formons leur produit $X\times X'\rightarrow Y\times Y'\leftarrow Z\times Z'$.

La suite spectrale d'Eilenberg-Moore associ\'ee au diagramme $X\times X'\rightarrow Y\times Y'\leftarrow Z\times Z'$ s'identifie au produit tensoriel des suites spectrales associ\'ees \`a $X\rightarrow Y\leftarrow Z$ et $X'\rightarrow Y'\leftarrow Z'$.
En particulier pour tout diagramme $X\rightarrow Y\leftarrow Z$, les diagonales $X\to X\times X$, $Y\to Y\times Y$, etc. font de la suite spectrale associ\'ee au diagramme une suite spectrale de $\F_p$-alg\`ebres bigradu\'ees compatible avec la structure multiplicative de l'aboutissement.

\begin{rem} Chaque espace profini formant le diagramme $X\rightarrow Y\leftarrow Z$ est la limite filtrante de ses quotients finis simpliciaux.
En utilisant le lemme \ref{reind} on obtient que le diagramme lui m\^eme est la limite filtrante de diagrammes $X(i)\rightarrow Y(i)\leftarrow Z(i)$ entre ensembles finis simpliciaux.
La suite spectrale d'Eilenberg-Moore en cohomologie {\it modulo} $p$ associ\'ee \`a $X\rightarrow Y\leftarrow Z$ s'identifie alors \`a la colimite (filtrante) des suites spectrales d'Eilenberg Moore classiques associ\'ees aux diagrammes $X(i)\rightarrow Y(i)\leftarrow Z(i)$ de fa\c{c}on compatible avec sa structure.
\end{rem}

\section{Rappels sur la cohomologie {\it modulo} $p$ des espaces fonctionnels de source le classifiant du groupe cyclique d'ordre $p$}
\label{theorie:lannes}

On note $\B\Z/p$ le classifiant du groupe cyclique d'ordre $p$. C'est un ensemble fini simplicial repr\'esentant le foncteur qui associe \`a un espace profini l'ensemble de ses $1$-cocycles continus \`a coefficients dans $\Z/p$.

Rappelons que $\u$ d\'esigne la cat\'egorie des modules instables sur l'alg\`ebre de Steenrod {\it modulo} $p$ et notons $\H$ la cohomologie {\it modulo} $p$ de $\B\Z/p$.
Le foncteur $\u\to\u$, $M\mapsto \H\otimes M$ admet un adjoint \`a gauche
not\'e $\T$ (voir \cite{La}).
Pour tout espace profini $X$, l'\'evaluation $\B\Z/p\times\hom(\B\Z/p,X)\to X$ induit un morphisme naturel $\T\H^{*}X\to \H^{*}\hom(\B\Z/p,X)$.

\begin{theo}[\cite{La}, \cite{DrS}] Pour tout espace $p$-fini $X$ le morphisme $\T\H^{*}X\to \H^{*}\hom(\B\Z/p,X)$ est un isomorphisme.
\end{theo}

Morel en d\'eduit le th\'eor\`eme suivant :

\begin{theo}{\rm\cite{Mo}}\qua Soit $X$ un espace profini fibrant ; alors le morphisme $\T\H^{*}X\to \H^{*}\hom(\B\Z/p,X)$ est un isomorphisme.
\end{theo}

Indiquons comment le dernier th\'eor\`eme se d\'eduit du pr\'ec\'edent :
on observe d'abord que lorsque $X$ est fibrant dans $\wS$ le morphisme de $X$ dans sa r\'esolution fibrante $\R X=\lim \widehat{X}(-)$ est une \'equivalence d'homotopie simpliciale.
On en d\'eduit que pour $X$ fibrant le morphisme $\hom(\B\Z/p,X)\to\hom(\B\Z/p,\R X)$ est une \'equivalence d'homotopie donc induit un isomorphisme en cohomologie {\it modulo} $p$ continue.
On observe ensuite que le morphisme $\T\H^{*}\R X\to \H^{*}\hom(\B\Z/p,\R X)$ est la colimite des morphismes $\T\H^{*}\widehat{X}(i)\to \H^{*}\hom(\B\Z/p,\widehat{X}(i))$ o\`u les $\widehat{X}(i)$ sont des espaces $p$-finis.

Soit maintenant $\bar{\mathrm{H}}$ la cohomologie {\it modulo} $p$ r\'eduite de $\B\Z/p$. Le foncteur $\u\to\u$, $M\mapsto {\bar{\mathrm{H}}}\otimes
M$ admet \'egalement un adjoint \`a gauche, le foncteur de Lannes \emph{r\'eduit} not\'e $\bar{\mathrm{T}}$.
La d\'ecomposition $\H = \fp \oplus \bar{\mathrm{H}}$ induit pour toute module instable $M$ un scindement canonique~:
$$\mathrm{T}M=M \oplus \bar{\mathrm{T}}M\  .$$
On peut comme pr\'ec\'edemment donner une interpr\'etation g\'eom\'etrique de $\bar{\T}$ :

Soit $X$ un espace profini, $\R X$ le remplacement fibrant de $X$ et $Y$ la cofibre du morphisme $\R X\to\hom(\B\Z/p,\R X)$ induit par le morphisme $\B\Z/p\to \pt$.

L'inclusion du point dans $\B\Z/p$ fait de $\R X$ un r\'etracte de $\hom(\B\Z/p,\R X)$ de sorte qu'on a une suite exacte courte
$$0\to \bar{\H}\r^*Y\to \H^*\hom(\B\Z/p,\R X)\to \H^*X\to 0\  .$$
L'isomorphisme $\T\H^*X\to \H^*\hom(\B\Z/p,\R X)$ induit alors un isomorphisme naturel
$$\bar{\T}\H^*X\to \bar{\H}\r^*Y\ .$$

\begin{cor}
\label{interpretation:geometrique}
Soit $M$ un module instable. Si $M$ est la cohomologie continue r\'eduite d'un espace profini $X$, il est de m\^eme pour $\bar{\mathrm{T}} M$.
\end{cor}

\section{Filtrations de la cat\'egorie $\u$
~(voir \cite{Sc1} et \cite{Ku})}
\label{filtrations}

Dans toute cette section, on se
sp\'ecialise au cas $p=2$.

Rappelons qu'une sous-cat\'egorie
ab\'elienne $\mathcal{D}$
d'une cat\'egorie ab\'elienne $\mathcal{C}$ est dite de Serre, si pour toute suite exacte courte
$A\lra B\lra C$
telle que $A$ et $C$ sont dans $\mathcal{D}$, alors $B$ est \'egalement dans $\mathcal{D}$.

\subsection{La filtration nilpotente}\label{filnil}

Soit $M$ un module instable et $m$ un
entier naturel. Pour tout \'el\'ement $x$ de
$M$, on note  $|x|$ le degr\'e de $x$. On d\'efinit des op\'erateurs
$$\Sq _m~: M \lra M, x \longmapsto {\Sq}^{|x|-m}x\  .$$
On dit qu'un \'el\'ement $x\in M$ est $s$-nilpotent si~:
$$\forall\, 0\leq t < s , \exists\, c \in \mathbb{N} , (\Sq _t)^c x=0 \  .$$
Dans le cas o\`u $x$ est $s$-nilpotent, mais n'est pas $(s+1)$-nilpotent, on dit
que $x$
est
exactement $s$-nilpotent.

Un \'el\'ement exactement $0$-nilpotent est dit r\'eduit. Un module instable est dit r\'eduit si tous ses \'el\'ements sont r\'eduits. En termes plus simples, un module instable $M$ est r\'eduit si et seulement si l'application
$$\Sq _0~: M \lra M, x \longmapsto {\Sq}^{|x|}x\  $$
est injective. Dans le cas o\`u $M$ est muni d'une structure d'alg\`ebre compatible, l'\'el\'evation au carr\'e co\"{\i}ncide avec $\Sq _0$ et $M$ est r\'eduit comme module instable si et seulement si $M$ est r\'eduit comme alg\`ebre, \ie n'a pas d'\'el\'ements nilpotents.

Un module instable est $s$-nilpotent si tous ses \'el\'ements sont au moins $s$-nilpotents.
La sous-cat\'egorie pleine de $\u$ dont les objets sont les modules $s$-nilpotents
est de Serre et stable par colimites.
Cette sous-cat\'egorie, not\'ee $\nil _s$, co\"{\i}ncide avec
la plus petite sous-cat\'egorie de Serre stable par colimite qui contient tous les modules qui
sont des suspensions $s$-i\`emes.

Soit $M$ un module instable. On d\'efinit
$M_s$ comme le plus grand sous-module de
$M$ qui est dans $\nil _s$. Les sous-modules $\{M_s\}_{s\in\n}$ forment une
filtration d\'ecroissante naturelle et
s\'epar\'ee de $M$.
Le sous-quotient $M_s /M_{s+1}$ est la suspension $s$-i\`eme ${\Sigma}^s {\R}_s M$ d'un module r\'eduit ${\R}_s M$.

Soient $M$ et $N$ deux modules
instables. On a $(M\otimes N)_n \cong
{\oplus}_{i+j=n} M_i \otimes N_j$. Par
cons\'equent, on a pour tout $n$ la
formule $R _n (M\otimes
N)\cong{\oplus}_{i+j=n} {\R}_ i M \otimes {\R}_j
N$ (voir \cite[Proposition 2.5, (1)]{Ku}).
En particulier on a un isomorphisme naturel $\R_n\Sigma M\cong\R_{n-1}M$.

On dispose \'egalement d'une
caract\'erisation du degr\'e de nilpotence d'un module
instable en terme de foncteur de Lannes:
un module instable $M$ est $s$-nilpotent si et
seulement si pour tout $n$,
${\mathrm{T}} ^n M$ est $(s-1)$-connexe ({\it cf} \cite[D\'efinition-Proposition 6.1.1, p. 139]{Sc1}).

\subsection{La filtration de Krull}

La filtration de Krull $\{\mathcal{U}_n\}_{n\in \mathbb{N}}$ de la
cat\'egorie $\u$ est d\'efinie par :
$\u_n$ est la sous-cat\'egorie pleine de $\u$ form\'ee des objets annul\'es par $\bar{\T}^n$.

C'est une filtration croissante
par des sous-cat\'egories de Serre
stables par colimites et suspensions,
qui est exhaustive au sens o\`u la plus
petite sous-cat\'egorie de $\u$ stable par
colimite et contenant tous les $\mathcal{U} _n$ pour tout entier $n$ est
$\u$ elle-m\^eme.

La filtration de Krull induit sur chaque module instable une filtration
croissante naturelle et compl\`ete par des sous-modules qui sont dans $\u _n$.
Dans le cas des modules r\'eduits, cette filtration co\"{\i}ncide avec la \emph{filtration par le
poids}, dont la d\'efinition suit.

\begin{defi}[~\em(Poids d'un module r\'eduit)]
\begin{enumerate}

\item Soit $n$ un entier naturel, on d\'efinit $\alpha (n)$ comme le nombre de $1$
dans l'\'ecriture diadique de $n$.

\item Un module instable r\'eduit est de poids inf\'erieur ou \'egal \`a $n$ s'il est
nul dans les degr\'es $i$ tels que $\alpha (i) >n$.

\end{enumerate}

\end{defi}

Un module r\'eduit est filtr\'e par ses sous-modules de poids $n$ maximaux.
Cette filtration est exhaustive. Le poids d'un module instable r\'eduit $M$ est not\'e
$w(M)$.
La caract\'erisation suivante est extraite de \cite{fs}.

\begin{pro}[Voir \cite{fs}]
\label{carac:un}
Un module r\'eduit est dans $\u _n$ si et seulement s'il est de poids inf\'erieur ou
\'egal \`a $n$.
\end{pro}

 Remarquons enfin que si $M$ est
dans $\u _n$ et si $N$ est dans $\u _m$, alors $M\otimes
N$ est dans $\u _{n+m}$.

\section{D\'emonstration du th\'eor\`eme \ref{theo:dg}}

Dans cette section, on fixe $p=2$.

On d\'emontre dans cette section le th\'eor\`eme \ref{theo:dg} en reprenant la preuve de L. Schwartz.
Ceci d\'emontre la
{\it conjecture de non r\'ealisation forte} pour $p=2$ en toute g\'en\'eralit\'e puisque la cohomologie {\it modulo} $2$ d'un espace est aussi la cohomologie {\it modulo} $2$ continue de son compl\'et\'e profini.

On suppose qu'il existe un espace
profini dont la cohomologie {\it modulo} $2$ continue n'est pas localement
finie, et qui est dans cran fini de la filtration de Krull. On va trouver une contradiction.
La premi\`ere \'etape dans ce sens est {\it la  r\'eduction de N. Kuhn}~:~on va montrer qu'il
existe alors un espace profini dont la
cohomologie continue est dans $\u _1$ mais n'est
pas dans $\u _0$ (\ie n'est pas
localement finie).

\subsection{La r\'eduction de N. Kuhn}

\begin{pro}[R\'eduction de Kuhn] S'il existe un entier $n\geq 1$ et un espace profini dont la cohomologie \emph{modulo} $2$ continue est dans $\u_n$ mais n'est
pas dans $\u_{n-1}$, alors il existe un espace profini dont la cohomologie {\it modulo} $2$ continue est
dans $\u_1$
mais n'est pas dans $\u_0$.
\end{pro}

{\bf D\'emonstration}\qua Soit $X$ un espace profini dont la cohomologie continue $M$ est dans ${\u}_{n}$ mais n'est pas dans ${\u}_{n-1}$ et supposons $n>1$.
Alors $\bar{\T}M$ est dans $\u_{n-1}$ mais pas dans $\u_{n-2}$, et est la cohomologie continue r\'eduite d'un espace profini d'apr\`es le corollaire \ref{interpretation:geometrique}.
On montre ainsi l'\'enonc\'e par une r\'ecurrence descendante.
\fin

\subsection{Construction des classes ${\alpha}_{i,d}$ (d'apr\`es \cite{Sc2})}

On s'est ramen\'e \`a consid\'erer un espace profini $X$ dont la cohomologie continue est dans $\u _1$
mais n'est pas localement finie. On utilise la structure tr\`es particuli\`ere de ces
modules pour construire des classes de la cohomologie continue d'un espace de lacets it\'er\'es de $X$
qui ont des propri\'et\'es contradictoires.

Soit $X$ un espace profini dont la cohomologie continue est dans $\u _1$ mais n'est pas
localement finie. Le module instable
$\bar{\mathrm{T}}\mathrm{H}^* X$ est non
nul par d\'efinition. Soit $(d-1)$ sa
connexit\'e.

Remarquons que si un
espace profini $X$ a sa cohomologie continue dans $\u _1$ mais pas
dans $\u _0$, il en est de m\^eme pour
toutes ses suspensions. En effet, le foncteur $\bar{\mathrm{T}}$ est exact et commute aux
suspensions et aux colimites.

On peut donc supposer que $X$
est une suspension et, par cons\'equent,
on peut supposer que $d\geq 1$ et que tous les
cup-produits sont nuls dans la
cohomologie continue r\'eduite de $X$.

On d\'efinit pour tout entier $n$ le $n$-squelette $\Sk_n Z$ d'un espace profini $Z$ comme le sous-espace profini engendr\'e par les simplexes non d\'eg\'en\'er\'es de dimension inf\'erieure ou \'egale \`a $n$.
Si $Z$ est la limite d'un diagramme filtrant d'ensembles finis simpliciaux $Z(-)$, alors $\Sk_n Z$ est la limite des $n$-squelettes des $Z(i)$.

On voit que la cohomologie continue r\'eduite du quotient de $X$ par son $(d-1)$-squelette est un module instable $(d-1)$-connexe qui ne diff\`ere de la cohomologie continue de $X$ que par un module born\'e, donc localement fini. Par cons\'equent $\bar{\mathrm{T}} \mathrm{H}^* X=
\bar{\mathrm{T}} \bar{\mathrm{H}}^*(X/\Sk_{d-1} X)$. Ainsi $X/\Sk_{d-1} X$ est \'egalement un espace profini dont la cohomologie continue est dans $\u _1$ mais n'est pas localement finie. De plus $X/\Sk_{d-1} X$ est encore une suspension donc tous les cup-produits sont nuls dans sa cohomologie continue r\'eduite.

On peut donc aussi supposer que $\bar{\H}^* X$ est $(d-1)$-connexe. L'int\'er\^et de cette d\'emarche est la simplification suivante~: un module instable est $s$-nilpotent si et seulement si ${\mathrm{T}}^{n} M$ est $(s-1)$-connexe pour tout $n$ (voir la section \ref{filtrations}). En vertu du scindement ${\mathrm{T}} M =M \oplus \bar{\mathrm{T}} M$, il appara\^{\i}t que $M$ est $s$-nilpotent si et seulement si $M$ est $(s-1)$-connexe et $\bar{\mathrm{T}} ^{n} M$ est $(s-1)$-connexe pour tout $n$.
En particulier si $M$ est dans $\u_1$, $M$ est $s$-nilpotent si et seulement si $M$ et $\bar{\T}M$ sont $(s-1)$-connexes.

Il r\'esulte de toutes les consid\'erations pr\'ec\'edentes qu'on peut supposer que $\bar{\mathrm{H}} ^{*} X$ est dans $\u_1$ mais n'est pas localement finie, est $(d-1)$-connexe, est $d$-nilpotente et poss\`ede une structure d'alg\`ebre triviale.

Soit ${\Fun}$ le sous-module instable de $\mathrm{H}^*\B{\mathbb{Z}}/2$ engendr\'e par la classe de degr\'e un. La construction des classes $\alpha _{i,d}$ annonc\'ee dans l'introduction repose la proposition suivante.

\begin{pro}\label{u1}
Soit $M$ un module qui est dans $\u _1$ mais n'est pas  localement fini.
Soit $\eta$ l'unit\'e d'adjonction $M\to \bar{\T}M\otimes \bar{\H}\r^*\B\Z/2$.
Alors $\eta$ factorise par le sous-module $\bar{\mathrm{T}}  M \otimes {\Fun}$. De plus le noyau et le conoyau de $$\eta:M\to \bar{\T}M\otimes {\Fun}$$ sont localement finis.
\end{pro}

{\bf D\'emonstration}\qua Dans \cite{Sc2}, ce r\'esultat est d\'emontr\'e sous
l'hypoth\`ese que $M$ est de type fini.
Tout module instable est colimite filtrante de ses
sous-modules de type fini. Etant donn\'e que les
colimites filtrantes sont exactes, commutent aux
produits tensoriels ainsi qu'au foncteur $\bar{\mathrm{T}}$,
et du fait que $\u _0$ est stable par colimites,
il vient que ce r\'esultat est vrai en l'absence de
toute hypoth\`ese de finitude.
\fin

\begin{lem}
\label{const:classes}
Il existe un entier $i_d$ tel que pour tout $i\geq i_d$, il existe des classes
$\alpha_{i,d}$ dans $\bar{\mathrm{H}}^*
X$ qui v\'erifient~:

\begin{itemize}

\item{la classe $\alpha_{i,d}$ est de degr\'e $2^i+d \ ,$}

\item{la classe $\alpha_{i,d}$ est exactement $d$-nilpotente,}

\item{$\Sq ^{2^i} \alpha_{i,d}=  \alpha_{i+1,d}\  .$}

\end{itemize}

\end{lem}

{\bf D\'emonstration}\qua Posons $M=\bar{\H}\r^* X$. On applique la proposition
\ref{u1}. et on obtient une suite
exacte~:
$$0\lra L \lra M \stackrel{\eta}{\lra}
\bar{\mathrm{T}}M \otimes {\Fun}
\lra L' \lra 0$$
o\`u les modules instables $L$ et $L'$ sont localement finis.
Comme
$\bar{\mathrm{T}}M$ est $(d-1)$-connexe, il
existe $\alpha\in\bar{\mathrm{T}}M$ non nul et de degr\'e
$d$.

Etant donn\'e que $\bar{\mathrm{T}} M$ est localement fini, il
existe une borne enti\`ere $h$ telle que toute op\'eration
de Steenrod de degr\'e plus grand que $h$ annule
l'\'el\'ement $\alpha$ de $\bar{\mathrm{T}} M$, \ie
${\mathcal{A}}_{2}^{\geq h}.\alpha =0$.

Par cons\'equent, pour tout
entier $i$ tel que $ 2^i \geq 2h$, on a
par la formule de Cartan,
\begin{eqnarray}
\label{action1}
{\Sq}^{2^i} ~\alpha \otimes u^{2^i} =
\alpha \otimes u^{2^{i+1}} \ .
\end{eqnarray}
Soit $\kappa$ un entier tel que $
2^{\kappa} \geq 2h$.

Il suit de (\ref{action1}) et du fait que $L'$
est localement fini qu'il existe une
borne $\kappa ' \geq \kappa$, telle que les classes
$\alpha\otimes {u}^{2^i}$ de $\bar{\mathrm{T}}M \otimes {\Fun}$
sont d'image nulle dans $L'$ pour
$i>\kappa '$. Ainsi, l'\'el\'ement $\alpha \otimes
{u}^{2^i}$ provient pour $i\geq\kappa '$ d'un \'el\'ement ${\alpha}_{i,d}'$
de $M$, \ie $\eta ({\alpha}_{i,d}')=\alpha\otimes {u}^{2^i}\ $.

On choisit alors~:
\begin{enumerate}

\item $i_d = \kappa '$,

\item ${\alpha}_{{i_d},d}={\alpha
'}_{{i_d},d}~$,

\item pour tout $i\geq i_d$,
${\alpha}_{i+1,d}={\Sq}^{2^i}~{\alpha}_{i,d}~$.

\end{enumerate}

Ces classes v\'erifient le premier et le troisi\`eme
point du lemme \ref{const:classes} par construction. Le second point d\'ecoule du fait
suivant : La formule de Cartan, l'action des carr\'es de Steenrod ${\Sq}^i$ sur ${\Fun}$ et le fait que $2^{\kappa '} \geq 2^{\kappa}\geq 2h$ entra\^{\i}nent
$$\forall\,c\geq 0,\forall\,t<d,\eta ({\Sq}_{t}^{c}~{\alpha}_{i,d})=
{\Sq}_{t}^{c}~({\alpha}\otimes {u}^{2^i})=
({\Sq}_{t}^{c}~{\alpha})\otimes {u}^{2^i}\ .$$
Or le degr\'e de
${\Sq}_{t}^{c}~{\alpha}$ est une fonction
strictement croissante de $c$ car $t<d$.
Donc, pour $c$ suffisamment grand, $\eta
({\Sq}_{t}^{c}~{\alpha}_{i,d})=0$ car
$\bar{\mathrm{T}} M$ est
localement fini. De ce fait
${\Sq}_{t}^{c}~{\alpha}_{i,d}$ provient de $L$ qui est
localement fini.
Donc \`a nouveau $\Sq_t^c\,\alpha_{i,d}$ est nul pour $c$ suffisamment grand, c'est-\`a-dire ${\alpha}_{i,d}$ est $t$-nilpotent.

Pour $t=d$, on a
$$\eta({\Sq}_{d}~{\alpha}_{i,d})=
{\Sq}_{d}~({\alpha}\otimes {u}^{2^i})
= {\alpha}\otimes {u}^{2^{i+1}}\neq 0\  $$
et donc pour tout entier $c$,
$$\eta({\Sq}_{d}^{c}~{\alpha}_{i,d})=
{\Sq}_{d}^{c}~({\alpha}\otimes {u}^{2^{i}})
= {\alpha}\otimes {u}^{2^{i+c}}\neq
0\  $$
ce qui montre que ${\Sq}_{d}^{c}~{\alpha}_{i,d}$ est
non nul quelque soit $c$.
On a montr\'e
que ${\alpha}_{i,d}$ est exactement $d$-nilpotente, ce
qui d\'emontre le deuxi\`eme point du lemme.
\fin

\subsection{Construction des classes ${\alpha}_{i,\ell}$}
\label{def:classes}

Pour $0\leq \ell \leq d$, et pour tout
$i\geq i_d$, on d\'efinit des classes
${\alpha}_{i,\ell}$ dans la cohomologie de ${\Omega}^{d-
\ell } X$ de la mani\`ere suivante~:

\begin{itemize}

\item les classes ${\alpha}_{i,d}$ sont d\'efinies par
le lemme \ref{const:classes},

\item la classe ${\alpha}_{i,\ell -1}$ est la
d\'esuspension de l'image de ${\alpha}_{i,\ell}$ par
l'homorphisme de coin $\bar{\H}^{*}{\Omega}^{d-\ell }X\lra \Sigma
\bar{\H}^{*}{\Omega}^{d-\ell +1}X\ $.
\end{itemize}

\begin{lem}
\label{prop:classes}
Les classes $\alpha_{i,\ell}$ construites pr\'ec\'edemment ont les propri\'et\'es
suivantes, pour $d \geq \ell \geq
0$ et pour tout $i\geq i_d$~:

\begin{enumerate}

\item{la classe $\alpha_{i,\ell}$ est de degr\'e $2^i+\ell$,}

\item{la classe $\alpha_{i,\ell}$ est exactement $\ell$-nilpotente,}

\item{$\Sq ^{2^i} \alpha_{i,\ell}=  \alpha_{i+1,\ell}$,}

\item\label{point4}{il existe un entier $i_{\ell}$ tel que
pour tout $i\geq i_{\ell}$, le cup-carr\'e
$\alpha_{i,\ell} \cup \alpha_{i,\ell}$
est nul.}

\end{enumerate}

\end{lem}

{\bf D\'emonstration}\qua Les trois premi\`eres affirmations du lemme sont
faciles~ :

\begin{itemize}

\item le premier point r\'esulte des d\'efinitions,

\smallskip\item le troisi\`eme point est cons\'equence des
propri\'et\'es de compatibilit\'e de la suite spectrale
d'Eilenberg-Moore aux op\'erations de Steenrod,

\smallskip\item le deuxi\`eme point est cons\'equence du fait qu'on a des
monomorphismes (voir la proposition  \ref{mono}) qui sont induits par le morphisme de coin it\'er\'e~:
$${\R}_{d} \bar{\mathrm{H}}^{*}X
\hra {\R}_{{d-1}} \bar{\mathrm{H}}^{*}\Omega X
\hra \ldots  \hra {\R}_{\ell} \bar{\mathrm{H}}^{*}{\Omega}^{d- \ell} X  \hra \ldots
\hra {\R}_{0} \bar{\mathrm{H}}^{*}{\Omega}^{d}X$$
En effet, soit $\bar{\alpha} _{i,\ell}$ l'image de ${\alpha _{i,\ell}}$ dans ${\R}_{\ell}
\bar{\mathrm{H}} ^{*}{\Omega}^{d-\ell }
X$. Les classes $\bar{\alpha} _{i,\ell}$ sont les images successives
de $\bar{\alpha} _{i,d}$ par les monomorphismes pr\'ec\'edents, et donc $\alpha _{i,\ell}$ est non nulle pour
tout $\ell\leq d$. Une r\'ecurrence sur la proposition \ref{filtration:nil:lacet:un} montre
d'autre part que $\bar{\mathrm{H}} ^{*}{\Omega}^{d-\ell } X$ est au moins
$\ell$-nilpotent. Par cons\'equent les classes $\alpha _{i,\ell}$ sont au moins
$\ell$-nilpotentes et r\'eduisent non trivialement dans  ${\R}_{\ell}
\bar{\mathrm{H}} ^{*}{\Omega}^{d-\ell } X$. C'est dire qu'elles
sont exactement $\ell$-nilpotentes.

\end{itemize}

Le point \ref{point4} se d\'emontre par une r\'ecurrence
descendante, que l'on renvoie au prochain paragraphe.
\fin

Si l'on admet ce quatri\`eme point, la d\'emonstration de la conjecture de N. Kuhn
est termin\'ee, car les deuxi\`eme et quatri\`eme points du lemme  \ref{prop:classes} sont
contradictoires pour $l=0$.

\subsection{D\'emonstration du point \ref{point4}  du lemme \ref{prop:classes}}
\label{cupnul}

Remarquons tout d'abord que par application it\'er\'ee des propositions \ref{filtration:nil:lacet:zero} et \ref{filtration:nil:lacet:un} et du corollaire \ref{corwin}, on a le lemme suivant~:

\begin{lem}
\label{thelem}
Pour $0 \leq \ell \leq d$, on a~:
\begin{itemize}

\item
Le module instable $\bar{\H}^{*}{\Omega}^{d-\ell} X$ est (au moins)
$\ell$-nilpotent,

\item le module instable ${\R}_i {\F}_{-1} \bar{\mathrm{H}}^{*} {\Omega}^{d-\ell } X$ est dans $\u _1$ pour
$i\leq 2\ell\ ,$

\item  le module instable ${\R}_{i} {\F}_{-2}
\bar{\mathrm{H}} ^{*}{\Omega}^{d-\ell }
X$ est dans $\u _2$ pour
$i\leq 2\ell$,

\item  si $\ell \geq 1$ et si $i\leq 2\ell -1$, le module ${\R}_{i} {\F}_{-2} \bar{\mathrm{H}} ^{*}{\Omega}^{d-\ell }
X$ est en fait dans $\u _1$.

\end{itemize}
\end{lem}

Pour d\'emontrer le point
\ref{point4}  du lemme
\ref{prop:classes}, on proc\`ede par r\'ecurrence descendante sur $\ell$.

Soit $\ell\geq 1$. Supposons que pour $i\geq i_\ell $,
${\alpha}_{i,\ell }\cup{\alpha}_{i,\ell }=0$.

On consid\`ere la suite spectrale
d'Eilenberg-Moore qui relie la
cohomologie de ${\Omega}^{d-\ell} X$ \`a
celle de ${\Omega}^{d-\ell +1}X $.

Alors ${\alpha}_{i,\ell }\otimes{\alpha}_{i,\ell }$ d\'efinit
un $1$-cycle, qui est un cycle permanent. En
effet, ${\alpha}_{i,\ell}\otimes{\alpha}_{i,\ell}$ est
exactement $2\ell$-nilpotente et l'image des
diff\'erentielles est constitu\'ee
d'\'el\'ements de degr\'e de nilpotence
strictement sup\'erieur \`a $2\ell$. Soit ${\omega}_{i, \ell -
1}$ une classe de 
$\mathrm{H}^{*}{\Omega}^{d-\ell +1}X$ qui est d\'etect\'ee
par le cycle permanent induit par ${\alpha}_{i,\ell}\otimes{\alpha}_{i,\ell}$.

\begin{lem}
\label{nilp:omega}
La classe ${\omega}_{i, \ell -1}$ est au moins $(2\ell-2)$-nilpotente.
\end{lem}

{\bf D\'emonstration}\qua
Remarquons que cette condition est
trivialement v\'erifi\'ee si $\ell =1$.

Pour $\ell> 1$, la classe
${\omega}_{i, \ell -1}$ est au moins
$(\ell -1)$-nilpotente car le module
${\F}_{-2}\bar{\mathrm{H}}^{*} {\Omega}^{d-\ell +1}$ est $(\ell -1)$-nilpotent.

Soit $\displaystyle{\overline{{\omega}_{i,
\ell -1}}^s}$ l'image
de ${\omega}_{i, \ell -1}$ dans
${\R}_s \bar{\mathrm{H}}^{*} {\Omega}^{d-\ell +1} X$ pour $s\leq 2\ell -3$.
La classe $\displaystyle{\overline{{\omega}_{i, \ell -1}}^s}$ est de degr\'e $2^{i+1}+2\ell -2-
s$. On a
$$\alpha(2^{i+1}+2\ell -2-s)\geq 2\ ,$$

or ${\R}_s
\bar{\mathrm{H}}^{*} {\F}_{-2}{\Omega}^{d-\ell +1} X$ est dans $\u _1$, pour $\ell -
1\leq s \leq 2\ell
-3$ d'apr\`es le lemme \ref{thelem}. De ce fait, pour $\ell -1\leq s \leq 2\ell -3$,
$\displaystyle{\overline{{\omega}_{i, \ell
-1}}^s}$ est nulle et donc ${\omega}_{i, \ell -1}$
est au moins  $(2\ell-2)$-nilpotente.
\fin

\begin{lem}
\label{steen:omega}
On a l'\'egalit\'e
$${\Sq}^{2^i}~{\omega}_{i, \ell -1}={\alpha}_{i, \ell -1}\cup
{\alpha}_{i+1,\ell -1}\quad {\rm
mod~}({\F}_{-2} {\bar{\mathrm{H}}^* }{\Omega}^{d-\ell +1}
X)_{2\ell -2}\ .$$
\end{lem}
{\bf D\'emonstration}\qua
On a par la formule de Cartan, pour $i\geq i_{\ell}$
{\small
$${\Sq}^{2^i}~({\alpha}_{i,\ell}\otimes{\alpha}_{i,\ell})=
({\Sq}^{2^i}~{\alpha}_{i,\ell})\otimes{\alpha}_{i,\ell}
+{\alpha}_{i,\ell}\otimes({\Sq}^{2^i}~{\alpha}_{i,\ell})=
{\alpha}_{i+1,\ell}\otimes{\alpha}_{i,\ell}
+{\alpha}_{i,\ell}\otimes{\alpha}_{i+1,\ell}\ .$$}
Or ${\alpha}_{i+1,\ell}\otimes{\alpha}_{i,\ell}
+{\alpha}_{i,\ell}\otimes{\alpha}_{i+1,\ell}$ est un cycle
permanent qui n'est jamais l'image d'une diff\'erentielle
(pour des raisons de nilpotence, comme
pr\'ec\'edemment) et qui
d\'etecte
${\alpha}_{i+1,\ell -1}\cup{\alpha}_{i,\ell -1}$.

Ceci montre que la diff\'erence
$z={\Sq}^{2^i}~{\omega}_{i, \ell -1}-{\alpha}_{i, \ell -1}\cup
{\alpha}_{i+1,\ell -1}$ vit dans ${\F}_{-1}
\bar{\mathrm{H}}^{*}{\Omega}^{d-\ell +1} X$.

On conclut en remarquant que
ces classes sont de degr\'e $2^i +
2^{i+1}+ 2\ell -2$. Or pour $s\leq 2\ell
-2$, on a $\alpha (2^i +
2^{i+1}+ 2\ell -2- s) \geq 2$. Comme le
module ${\R}_{s} {\F}_{-1}
\bar{\mathrm{H}} ^{*}{\Omega}^{d-\ell+1 }
X$ est dans $\u _1$ pour $s\leq 2\ell -2$ d'apr\`es le lemme
\ref{thelem}, on en d\'eduit l'\'egalit\'e
annonc\'ee.
\fin

On d\'eduit du lemme
\ref{steen:omega} qu'on a la relation
$${\Sq}^{2^i}{\Sq}^{2^i}~{\omega}_{i, \ell
-1}={\alpha}_{i+1, \ell -1}\cup
{\alpha}_{i+1,\ell -1}\quad {\rm
mod~}({\F}_{-2} {\bar{\mathrm{H}}^* }{\Omega}^{d-\ell +1}
X)_{2\ell -2}\ .$$
Ceci r\'esulte de la formule de
Cartan, et du fait que ${\R}_{2\ell
-2}{\F}_{-2} {\bar{\mathrm{H}}^*
}{\Omega}^{d-\ell +1}$ est dans $\u_2$, par
le lemme \ref{thelem}.

On remarque alors que
\begin{enumerate}

\item la relation ${\Sq}^{2^i}{\Sq}^{2^i}{\omega}_{i,
\ell -1}={\alpha}_{i+1, \ell -1}\cup {\alpha}_{i+1,\ell -
1}\quad$ est vraie {\it modulo} des termes de degr\'e
de nilpotence strictement plus grand que $2\ell-2$,

\item l'image
$\displaystyle{\overline{{\Sq}^{2^i}{\Sq}^{2^i}{\omega}_{i, \ell -1}}}$ de
${\Sq}^{2^i}{\Sq}^{2^i}{\omega}_{i, \ell -1}$ dans
${\R}_{2\ell-2}{\F}_{-2} \bar{\mathrm{H}}^{*}{\Omega}^{d-\ell +1} X$ est
nulle.

\end{enumerate}

Le deuxi\`eme point provient du
r\'esultat suivant \cite[lemme 5.7, p. 554]{Sc3}~:

\begin{lem}
Pour tout entier $n>0$
$${\Sq}^{2^n}{\Sq}^{2^n}\in \bar{\mathcal{A}}(n){\Sq}^{2^n}\bar{\mathcal{A}}(n)$$
o\`u ${\mathcal{A}}(n)$ est la sous-alg\`ebre engendr\'ee par $\{ \mathrm{\Sq}^{2^i}, 0\leq i\leq n-1\}$,
et $\bar{\mathcal{A}}(n)$ est l'id\'eal des \'el\'ements strictement positifs de ${\mathcal{A}}(n)$.
\end{lem}

Cette d\'ecomposition de l'op\'eration de Steenrod
${\Sq}^{2^i}{\Sq}^{2^i}$ montre que\\
$\displaystyle{\overline{{\Sq}^{2^i}{\Sq}^{2^i}{\omega}_{i, \ell -1}}}$ est
de poids sup\'erieur \`a $3$, alors que ${\R}_{2\ell-
2}\mathrm{H}^{*}{\Omega}^{d-\ell +1} X$ est
dans $\u _2$, d'apr\`es le lemme
\ref{thelem}.

Ceci montre que ${\alpha}_{i+1, \ell -1}\cup
{\alpha}_{i+1,\ell -1}$ est au moins $(2\ell -1)$-nilpotent et donc qu'il existe un
entier $c$ tel que
$$0={\Sq}_{2\ell -2}^{c}({\alpha}_{i +1, \ell -1}\cup
{\alpha}_{i+1,\ell -1})={\alpha}_{i+c
+1, \ell -1}\cup
{\alpha}_{i+c+1,\ell -1}\  .$$

On peut donc choisir $i_{\ell
+1}=i_{\ell} +c+1$, et on aura
pour tout $i\geq i_{\ell +1}$, ${\alpha}_{i, \ell -1}\cup
{\alpha}_{i,\ell -1}=0\  .$

\appendix

\section{Appendice}

\subsection{Filtration de Krull et filtration nilpotente}

Rappelons qu'on note pour $M$ un module instable et $s$ un entier $M_s$ le sous-module instable de $M$ form\'e des \'el\'ements (au moins) $s$-nilpotents et $\R_sM$ le module $\Sigma^{-s}(M_s/M_{s+1})$ (voir la section \ref{filnil}).

La proposition suivante affirme {\it grosso
modo} que les foncteurs ${\R}_{s}$
sont exacts \`a gauche {\it modulo} $\nil$.

\begin{pro}
\label{exac:gauche}
Soit
$$0\lra A\lra B\lra C \lra 0$$
une suite exacte dans $\u$ et soit $L$ le conoyau de
${\R}_{s} A\lra {\R}_{s}
B$.
Pour tout entier $s$, le morphisme ${\R}_{s} A\lra {\R}_{s} B$
est un monomorphisme et le noyau de $L \lra {\R}_{s} C $
est constitu\'e d'\'el\'ements nilpotents.
\end{pro}

{\bf D\'emonstration}\qua
On chasse dans le diagramme commutatif suivant, dont les trois colonnes sont exactes.
\begin{eqnarray}\label{nilpot}
\begin{array}{ccccc} A_{s+1} & \longrightarrow & B_{s+1} & \longrightarrow & C_{s+1} \\
\downarrow & & \downarrow & & \downarrow \\
A_s & \longrightarrow & B_s & \longrightarrow & C_s \\
\downarrow & & \downarrow & & \downarrow \\
\Sigma^s {\R}_s A & \longrightarrow & \Sigma^s {\R}_s B & \longrightarrow & \Sigma^s {\R}_s C
\end{array}
\end{eqnarray}


{\it  Le morphisme ${\R}_{s} A\lra {\R}_{s}B$ est un monomorphisme~:}
Soit $x$ un \'el\'ement de $\ker ({\R}_{s} A\lra {\R}_{s}
B)$. Soit $x_1$ un \'el\'ement de ${A}_{s}$ qui se
projette sur $\Sigma ^s x\in \Sigma ^s {\R}_{s} A $.

Soit $x_2$ l'image de $x_1$ dans ${B}_s$. L'\'el\'ement $x_2$ se projette sur $0$ dans $\Sigma ^s {\R}_{s} B $.

Par exactitude de la seconde colonne, l'\'el\'ement $x_2$ est dans $B_{s+1}$, \ie
$x_2$ est $(s+1)$-nilpotent donc pour un certain
entier $c$, on a $({\Sq}_{s})^{c}x_2 =0$.

Donc par injectivit\'e de
${A}_{s}\lra {B}_{s}$, il vient que
$({\Sq}_{s})^{c}x_1 =0$, \ie l'\'el\'ement $x_1$ est dans
$A_{s+1}$ soit \`a dire que $x_1$ est $(s+1)$-nilpotent .

 Par cons\'equent $x$ est nul, ce qui montre que
${\R}_{s}A \lra {\R}_{s}B$ est un
monomorphisme.

\smallskip
{\it Le noyau de $L \lra {\R}_{s} C $ est nilpotent~:}
Soit $x$ un \'el\'ement de $\ker({\R}_{s}B \lra {\R}_{s}C)$ et soit $x_1$ un rel\`evement de $\Sigma ^s x$ dans
${B}_s$. Soit $x_2$ l'image de $x_1$ dans ${C}_s$.

L'\'el\'ement $x_2$ se projette sur $0$ dans $\Sigma ^s
{\R}_{s} C$, donc $x_2$ est dans $C_{s+1}$, \ie $x_2$ est $(s+1)$-nilpotent et pour un certain
entier $c$, on a $({\Sq}_{s})^{c}x_2 =0$.

Par cons\'equent, l'image de
$({\Sq}_{s})^{c}x_1$ dans $C$ est nulle et par exactitude de la suite $A \lra B \lra C $,
il existe un \'el\'ement $x_3$ de $A$ dont l'image est
$({\Sq}_{s})^{c}x_1$.
Comme $({\Sq}_{s})^{c}x_1$ est au moins $s$-nilpotent et comme $A\to B$ est un monomorphisme, l'\'el\'ement $x_3$ est au moins $s$-nilpotent.

Soit ${\Sigma}^s x_4$ l'image de  $x_3$ dans
${\Sigma}^s {\R}_{s} A$. L'image de ${\Sigma}^s x_4$ dans ${\Sigma}^s
{\R}_{s} B$ est
$({\Sq}_{s})^{c}{\Sigma}^s x ={\Sigma}^s ({\Sq}_{0})^{c}x$.
Donc l'image de $x$ dans $L$ est nilpotent.
\fin

\begin{cor} \label{exac:gauche2}
Soient $0\to A\to B\to C\to 0$ une suite exacte dans $\u$.
On suppose que $A$ est $\ell$-nilpotent pour un entier $\ell\geq 1$ ; alors le morphisme $\R_s B\lra \R_s C$ est un isomorphisme pour tout entier $s<\ell$ (et un \'epimorphisme pour $s=\ell$).
\end{cor}

{\bf D\'emonstration}\qua
Le module $\R_s A$ est nul pour tout $s<\ell$. La proposition affirme alors que le noyau du morphisme $\R_s B\lra \R_s C$ est nilpotent pour tout $s<\ell$.
Or $\R_s B$ est un module r\'eduit, donc ce noyau est nul.

Pour achever la preuve, il suffit de montrer que le morphisme $B_s\to C_s$ est un \'epimorphisme pour tout $s\leq\ell$.
Soit donc $s$ un tel entier et $x$ un \'el\'ement $s$-nilpotent de $C$.
Il provient d'un \'el\'ement $x_1$ de $B$ et pour tout $s'<s$ il existe un certain entier $c$ et un \'el\'ement $x_2$ de $A$ tel que $(\Sq_{s'})^c x_1$ est l'image de $x_2$ par le morphisme $A\to B$.
Or $A$ est $\ell$-nilpotent donc $(\Sq_{s'})^{c'} x_2$ est nul pour $c'$ assez grand.
On en d\'eduit que $x_1$ est $s$-nilpotent, c'est-\`a-dire $x_1\in B_s$.
\fin

\begin{cor} \label{U-Nil}
Soient $0\to A\to B\to C\to 0$ une suite exacte de modules instables et $p,q,s$ trois entiers positifs. On suppose que $\R_sA$ est dans $\u_p$ et que $\R_sC$ est dans $\u_q$ ; alors $\R_sB$ est dans $\u_{\mathrm{max}\{p,q\}}$.
\end{cor}

{\bf D\'emonstration}\qua Le foncteur $\R_s$ transforme la suite exacte
$$0\to A\to B\to C\to 0$$ en une suite $$0\to\R_s A\to \R_s B\to \R_s C$$ exacte au sens de la proposition \ref{exac:gauche}.

Les modules $\R_s A$, $\R_s B$ et $\R_s C$ sont r\'eduits, ce qui nous autorise \`a appliquer la caract\'erisation \ref{carac:un} de la filtration de Krull en terme de filtration par le poids pour les modules r\'eduits. Il suffit donc de montrer que $\R_sB$ est de poids inf\'erieur ou \'egal \`a ${\mathrm{max}\{p,q\}}$, \ie que $\R_sB$ est nul dans les degr\'es $i$ tels que $\alpha (i) > {\mathrm{max}\{p,q\}}$. Rappelons que $\alpha (i)$ est le nombre de puissances de deux composant l'\'ecriture diadique de $i$.

Soit $x$ un \'el\'ement de $\R_sB$ en degr\'e $|x|$ tel que $\alpha (|x|) > {\mathrm{max}\{p,q\}}$.
L'image de $x$ par $\R_s B\to \R_s C$ est nulle car $\R_s C$ est dans $\u_q$ et donc de poids inf\'erieur ou \'egal \`a $q\leq {\mathrm{max}\{p,q\}}$. D'apr\`es la proposition \ref{exac:gauche}, il existe un entier $c$ tel que ${Sq}_{0}^{c} ~x$ est dans $\R_s A$.
Le degr\'e de ${Sq}_{0}^{c} ~x$ est $2^c |x|$ et $\alpha (2^c |x|)=\alpha (|x|)>p$. Or $\R_s A$ est dans $\u_p$ , donc de poids inf\'erieur ou \'egal \`a $p$. Il s'ensuit que ${Sq}_{0}^{c} ~x$ est nul, et donc $x$ est nul car $\R_s A$ est r\'eduit.
\fin

\subsection{Suite spectrale d'Eilenberg-Moore et filtration nilpotente}

Dans toute cette section, $X$
est un espace profini point\'e.
On va d\'egager
suivant \cite{Sc3} quelques propri\'et\'es
de compatibilit\'e entre la filtration
nilpotente de la cohomologie continue r\'eduite de $X$ et celle de son espace
de lacets $\Omega X$ en utilisant la
spectrale d'Eilenberg-Moore.

Le terme
${\E}_{r}^{-s,*}$ de cette suite
spectrale est
naturellement muni d'une structure de
module instable pour tous entiers $s$ et
$r\geq 1$. De plus, la diff\'erentielle
$\d_r$ est un morphisme de modules
instables
$$\d_r: {\E}_{r}^{-s-r,*} \lra
{\Sigma}^{r-1} {\E}_{r}^{-s,*}\  .$$
Enfin, la cohomologie continue r\'eduite de $\Omega X$ est
naturellement munie d'une filtration
croissante
{\small $$0= {\F}_0 \bar{\H}^* \Omega
X\subset {\F}_{-1} \bar{\H}^* \Omega X
\subset {\F}_{-2} \bar{\H}^* \Omega
X\subset \ldots \subset {\F}_{-s}
\bar{\H}^* \Omega X\subset \ldots
\subset \bar{\H}^* \Omega X$$}
par des sous-modules
instables telle qu'on a pour tout $s\geq 1$ un isomorphisme
$${\E}_{\infty}^{-s,*}\cong {\Sigma}^{s}({\F}_{-s}\bar{\H}^* \Omega X/{\F}_{-s+1}\bar{\H}^* \Omega X)\  .$$

Cette filtration est convergente :
$${\bigcup}_{i\in\n} ~{\F}_{-i}\bar{\H}^* \Omega X =\bar{\H}^* \Omega X \  .$$
(Voir la section \ref{EM}.)

\begin{pro}
\label{mono}
Soit $X$ un espace profini point\'e dont la cohomologie r\'eduite est dans $\nil_{\ell}$ pour un entier ${\ell}\geq 1$.
L'homomorphisme de coin induit pour $s\leq 2{\ell}-1$ un monomorphisme $$\R_s\bar{\H}^*X \hra \R_{s-1}\bar{\H}^*\Omega X\  .$$
\end{pro}

{\bf D\'emonstration}\qua
L'homomorphisme de coin
$\bar{\mathrm{H}}^{*}X\lra\Sigma
\bar{\mathrm{H}}^{*}\Omega X$ factorise de la mani\`ere
suivante:
$$\bar{\mathrm{H}}^{*}X\cong \E_1^{-1,*}
\era {\E}^{-1,*}_{\infty}\cong
\Sigma {\F}_{-1}\bar{\mathrm{H}}^{*}\Omega X
\hra \Sigma \bar{\mathrm{H}}^{*}\Omega X\  .$$

Le morphisme $\E_1^{-1,*}\lra\E_\infty^{-1,*}$ est la colimite sur $r$ des morphismes $\E_1^{-1,*}\lra\E_r^{-1,*}$, lesquels sont les compos\'es $\E_1^{-1,*}\to\E_2^{-1,*}\to\cdots\to\E_{r-1}^{-1,*}\to\E_r^{-1,*}$.

Chaque morphisme $\E_r^{-1,*}\lra\E_{r+1}^{-1,*}$ est le conoyau de la diff\'erentielle $\d_r:\Sigma^{1-r}\E_r^{-1-r,*}\lra\E_r^{-1,*}$.
Comme le module   
${\E}^{-1-r,*}_{r}$ est un sous-quotient de 
$({\bar{\mathrm{H}}^{*}X})^{\otimes (r+1)}$ qui est 
$(r+1)\ell$-nilpotent, le noyau du morphisme $\E_r^{-1,*}\lra\E_{r+1}^{-1,*}$ est form\'e d'\'el\'ements au moins $((r+1)\ell-(r-1))$-nilpotents donc au moins $2\ell$-nilpotents puisqu'on a $r\geq 1$ et $\ell\geq 1$.

Le corollaire \ref{exac:gauche2} montre alors que le morphisme $\R_s\E_r^{-1,*}\lra\R_s\E_{r+1}^{-1,*}$ est un isomorphisme pour tout $s<2\ell$.
Par cons\'equent, le morphisme $\R_s\E_1^{-1,*}\lra\R_s\E_{\infty}^{-1,*}$ est un isomorphisme pour $s<2\ell$.

On a donc un isomorphisme
${\R}_s \bar{\mathrm{H}}^* X \lra {\R}_{s-1} {\F}_{-1} \bar{\mathrm{H}}^* \Omega X$
qui se compose avec le monomorphisme
$\R_{s-1} {\F}_{-1} \bar{\mathrm{H}}^* \Omega X\lra {\R}_{s-1}\bar{\mathrm{H}}^{*}\Omega X$
(le foncteur ${\R}_{s-1}$ pr\'eserve les monomorphismes d'apr\`es la proposition
\ref{exac:gauche})
en le monomorphisme souhait\'e.
\fin

\begin{pro}
\label{filtration:nil:lacet:un}
Soient $X$ un espace profini point\'e et
${\ell}\geq 1$ un entier tels que~:

\begin{itemize}

\item la cohomologie r\'eduite de $X$ est ${\ell}$-nilpotente,

\item $\bar{\mathrm{H}}^* X / (\bar{\mathrm{H}}^* X)_{2{\ell}}$ est dans $\u _1$.

\end{itemize}

Alors~:

\begin{enumerate}

\item \label{un} la cohomologie r\'eduite de $\Omega X$ est $({\ell}-1)$-nilpotente,

\item \label{deux} pour $s\leq 2\ell -2$, le module
instable ${\R}_s {\F}_{-1} \bar{\mathrm{H}}^* \Omega
X$ est dans $\u
_1$,

\item \label{trois} pour $s\leq 2\ell -2$, le module
instable ${\R}_{s} {\F}_{-2} \bar{\mathrm{H}}^* \Omega X$ est dans $\u_2$.

\end{enumerate}
\end{pro}

{\bf D\'emonstration du point
\ref{un}}\qua
Pour tout entier $s$,
le terme ${\E}^{-s,*}_{\infty}$ est un
sous-quotient de $(\bar{\mathrm{H}}^{*}
X)^{\otimes s}$ qui est
$s\ell$-nilpotent. Donc pour tout entier
$s$, le module ${\F}_{-s}\bar{\mathrm{H}}^{*}\Omega X/{\F}_{-s+1}
\bar{\mathrm{H}}^{*}\Omega
X ={\Sigma}^{-s}{\E}^{-s,*}_{\infty}$ est
$(s\ell -s)$-nilpotent. On en d\'eduit
par r\'ecurrence, en utilisant que
$\nil_{\ell -1}$ est de Serre, que pour
tout entier $s$, le module instable ${\F}_{-s}
\bar{\mathrm{H}}^{*}\Omega X$ est
$(\ell -1)$-nilpotent. Par
suite, comme la suite spectrale est
convergente et $\nil_{\ell -1}$ stable
par colimite, le module
$\bar{\mathrm{H}}^{*}\Omega X$ est aussi
$(\ell -1)$-nilpotent.

{\bf D\'emonstration du point
\ref{deux}}\qua Le module
${\E}^{-1,*}_{\infty}=\Sigma {\F}_{-1}\bar{\mathrm{H}}^{*}\Omega X$
est un quotient de ${\mathrm{H}}^{*} X$ par des
\'el\'ements au moins $2\ell$-nilpotents. On a donc un \'epimorphisme
$${\E}^{-1,*}_{\infty}\era \bar{\mathrm{H}}^* X /
(\bar{\mathrm{H}}^* X)_{2{\ell}}$$
de noyau au moins $2\ell$-nilpotent. En vertu du corollaire
\ref{exac:gauche2}, on a pour $s\leq
2\ell -1$ un monomorphisme ${\R}_s
{\E}^{-1,*}_{\infty} \lra {\R}_s
(\bar{\mathrm{H}}^* X /
(\bar{\mathrm{H}}^* X)_{2{\ell}})$. Or  $\bar{\mathrm{H}}^* X /
(\bar{\mathrm{H}}^* X)_{2{\ell}}$ est
dans $\u _1$ par hypoth\`ese et donc aussi
${\R}_s (\bar{\mathrm{H}}^* X /
(\bar{\mathrm{H}}^* X)_{2{\ell}})$ (car $\u _1$ est de Serre). Il
suit que ${\R}_s
{\E}^{-1,*}_{\infty}$ est \'egalement dans
$\u _1$.

On obtient que
$${\R}_{s} {\E}^{-1,*}_{\infty}\cong
{\R}_{s-1} {\Sigma}^{-1}{\E}^{-1,*}_{\infty}
\cong {\R}_{s-1}{\F}_{-1}\bar{\mathrm{H}}^{*}\Omega X$$
est dans $\u _1$ pour tout $s\leq 2\ell
-2$, ce qui montre le second point de la
proposition.

{\bf D\'emonstration du point
\ref{trois}}\qua
La suite exacte
$$0\lra
{\F}_{-1}\bar{\mathrm{H}}^{*}\Omega X
\lra {\F}_{-2}\bar{\mathrm{H}}^{*}\Omega X
\lra {\Sigma}^{-2}{\E}^{-2,*}_{\infty}\lra 0
$$
et le corollaire \ref{U-Nil} assurent que pour montrer le troisi\`eme point du lemme, il suffit de montrer que pour tout $s\leq
2\ell -2$, les modules ${\R}_s {\F}_{-1}\bar{\mathrm{H}}^{*}\Omega X$ et ${\R}_s {\Sigma}^{-2}{\E}^{-2,*}_{\infty}$ sont dans $\u_2$.

Le fait que ${\R}_s {\F}_{-1}\bar{\mathrm{H}}^{*}\Omega X$ est dans $\u _2$ pour $s\leq
2\ell -2$ est cons\'equence du point \ref{deux} de cette proposition, car $\u _1$ est une sous-cat\'egorie de $\u_2$.

Il reste donc \`a montrer que ${\R}_s {\Sigma}^{-2}{\E}^{-2,*}_{\infty}$ est dans $\u_2$. Le module instable ${\E}^{-2,*}_{\infty}$ est un sous-quotient de
${\bar{\mathrm{H}}^{*}X}^{\otimes 2}$, \ie
${\E}^{-2,*}_{\infty}$ s'\'ecrit $C/B$, avec
$C\subset{\bar{\mathrm{H}}^{*}X}^{\otimes 2}$.
De plus, le module $B$ est au moins
$3\ell$-nilpotent pour la raison suivante. On a
$$ \d_r : {\E}^{-2-r,*}_{r} \lra{\Sigma}^{(r-1)}{\E}^{-2,*}_{r} \  ,$$
ce qui montre que $B$ est constitu\'e d'\'el\'ements de
degr\'e de nilpotence au moins \'egal \`a
$(r+2)\ell -(r-1)\geq 3\ell$. On en d\'eduit avec le corollaire \ref{exac:gauche2} que pour $s< 3\ell$, on a des monomorphismes:
$$ \R _s {\E}^{-2,*}_{\infty} \hra \R _s C \hra  {\R}_{s}({\bar{\mathrm{H}}^{*}X}^{\otimes 2}) \  .$$

Le module ${(\bar{\mathrm{H}}^{*}X)}^{\otimes 2}$ est $2\ell$-nilpotent par hypoth\`ese, et donc le plus petit $s$ tel que
$\R _s ({\bar{\mathrm{H}}^{*}X}^{\otimes 2})$ est non trivial est $s=2\ell$. Or ${\R}_{2\ell}({\bar{\mathrm{H}}^{*}X}^{\otimes
2})$ est isomorphe \`a ${({\R}_{\ell} \bar{\mathrm{H}}^{*}X)}^{\otimes
2}$ qui est dans $\u _2$ (voir la partie \ref{filtrations}). Donc ${\R}_{2\ell}{\E}^{-2,*}_{\infty}$ est dans $\u _2$.
Finalement, ${\R}_{2\ell -2} {\Sigma}^{-2}{\E}^{-
2,*}_{\infty} \cong  {\R}_{2\ell}{\E}^{-
2,*}_{\infty}$ est dans $\u _2$.
\fin

\begin{pro}
\label{filtration:nil:lacet:zero}
Soit $\ell >1$ un entier et soit $X$ un
espace profini point\'e. Si $\bar{\mathrm{H}}^* X$ est $\ell$-nilpotent et si $\bar{\mathrm{H}}^* X
/ (\bar{\mathrm{H}}^* X)_{2{\ell}}$ est dans $\u _1$, alors $\bar{\mathrm{H}}^*
\Omega X
/ (\bar{\mathrm{H}}^* \Omega X)_{2({\ell}-1)}$
est \'egalement dans $\u _1$.
\end{pro}

{\bf D\'emonstration}\qua
Le module instable ${\E}^{-s,*}_{\infty}$ est un
sous-quotient de ${\bar{\mathrm{H}}^{*}X}^{\otimes
s}$, donc ${\E}^{-s,*}_{\infty}$ est $s{\ell}$-nilpotent et
${\F}_{-s}{\bar{\mathrm{H}}^{*}\Omega X}/{\F}_{-s+1}{\bar{\mathrm{H}}^{*}\Omega X} =
{\Sigma}^{-s}{\E}^{-s,*}_{\infty}$ est $s({\ell}-1)$-nilpotent.
On montre par r\'ecurrence sur $s$ que
${\F}_{-s}{\bar{\mathrm{H}}^{*}\Omega X}/{\F}_{-1}{\bar{\mathrm{H}}^{*}\Omega X}$ est $2({\ell}-1)$-nilpotent.
Par convergence de la suite spectrale, le module ${\bar{\mathrm{H}}^{*}\Omega X}/{\F}_{-1}{\bar{\mathrm{H}}^{*}\Omega X}$ est aussi $2({\ell}-1)$-nilpotent.

On a une suite exacte:
$$0\lra
{\F}_{-1}\bar{\mathrm{H}}^{*}\Omega X
\lra \bar{\mathrm{H}}^{*}\Omega X
\lra \bar{\mathrm{H}}^{*}\Omega X/{\F}_{-1}\bar{\mathrm{H}}^{*}\Omega X \lra 0$$
dont le dernier terme $\bar{\mathrm{H}}^{*}\Omega
X/{\F}_{-1}\bar{\mathrm{H}}^{*}\Omega X$ est $2({\ell}-1)$-nilpotent. En particulier, pour $s\leq 2\ell -3$, le module $\R_s(\bar{\mathrm{H}}^{*}\Omega
X/{\F}_{-1}\bar{\mathrm{H}}^{*}\Omega X)$ est nul et donc dans $\u_0$.
De plus, d'apr\`es la proposition pr\'ec\'edente, pour tout $s\leq 2\ell -3$, le module ${\R}_s {\F}_{-1} \bar{\mathrm{H}}^* \Omega X$ est dans $\u _1$. 

Le corollaire \ref{U-Nil} assure que
pour tout $s\leq 2{\ell}-3$, le module ${\R}_{s}\bar{\mathrm{H}}^{*}\Omega X$
est dans $\u _1$. Comme $\u _1$ est
de Serre, on obtient que $\bar{\mathrm{H}}^{*}\Omega
X/{(\bar{\mathrm{H}}^{*}\Omega X)}_{2{\ell}-2}$ est dans $\u_1$.
\fin

\begin{cor}\label{corwin} Soit $\ell >1$ un entier et soit $X$ un
espace profini point\'e. Si $\bar{\mathrm{H}}^* X$ est $\ell$-nilpotent et si $\bar{\mathrm{H}}^* X
/ (\bar{\mathrm{H}}^* X)_{2{\ell}}$ est dans $\u _1$, alors pour $s\leq 2\ell -3$, le module ${\R}_s {\F}_{-2}\bar{\mathrm{H}}^* {\Omega} X$ est dans $\u_1$.
\end{cor}

{\bf D\'emonstration}\qua Le module ${\R}_s {\F}_{-2}\bar{\mathrm{H}}^* {\Omega} X$ est un sous-module de ${\R}_s \bar{\mathrm{H}}^* {\Omega} X$ d'apr\`es la proposition \ref{exac:gauche}. Or pour $s\leq 2\ell -3$, le module ${\R}_s \bar{\mathrm{H}}^* {\Omega} X$ est isomorphe \`a ${\R}_s (\bar{\mathrm{H}}^* {\Omega}X / (\bar{\mathrm{H}}^* {\Omega}X)_{2{\ell}-2})$ , d'apr\`es le corollaire \ref{exac:gauche2}. Or ${\R}_s (\bar{\mathrm{H}}^* {\Omega} X / (\bar{\mathrm{H}}^* {\Omega} X)_{2{\ell}-2})$ est dans $\u_1$ comme sous-quotient de $\bar{\mathrm{H}}^* {\Omega}X / (\bar{\mathrm{H}}^* {\Omega}X)_{2{\ell -2}}$ qui d'apr\`es la proposition \ref{filtration:nil:lacet:zero} est dans $\u_1$.
\fin

\bibliographystyle{gtart}

\footnotesize\renewcommand{\refname}{R\'ef\'erences}

\Addresses\recd

\end{document}